\documentclass[10pt,final]{siamltex}
\usepackage{amssymb,amsmath,amsfonts}
\usepackage{psfig}
\usepackage{index}

\newtheorem{remark}[theorem]{Remark}

\newtheorem{example}[theorem]{Example}

\newtheorem{observation}[theorem]{Observation}
\newtheorem{comment}[theorem]{Comment}
\newtheorem{question}[theorem]{Question}

\def\define{\stackrel{\tiny def}=}
\def\definemath{\stackrel{\text{\tiny \text{def}}}=}
\def\definemathit{\stackrel{\text{\tiny \text{{\em def}}}}=}
\def\notstronglesseq{\mbox{$\hspace{0.055in}\bigcirc\hspace{-0.117in}{\scriptstyle{\not\preccurlyeq}}\hspace{0.055in}$}}
\def\notweaklesseq{\mbox{$\hspace{0.055in}\bigcirc\hspace{-0.11in}{\scriptstyle{\not\leqslant}}\hspace{0.055in}$}}
\def\stronglesseq{\mbox{$\hspace{0.055in}\bigcirc\hspace{-0.117in}{\scriptstyle{\preccurlyeq}}\hspace{0.055in}$}}
\def\weaklesseq{\mbox{$\hspace{0.055in}\bigcirc\hspace{-0.11in}{\scriptstyle{\leqslant}}\hspace{0.055in}$}}
\def\oplus{\mbox{$\hspace{0.055in}\bigcirc\hspace{-0.122in}+\hspace{0.055in}$}}
\def\otimes{\mbox{$\hspace{0.055in}\bigcirc\hspace{-0.122in}\times\hspace{0.055in}$}}
\def\ozero{\mbox{$\hspace{0.055in}\bigcirc\hspace{-0.094in}{\scriptstyle{0}}\hspace{0.055in}$}}
\def\oone{\mbox{$\hspace{0.055in}\bigcirc\hspace{-0.095in}{\scriptstyle{1}}\hspace{0.055in}$}}

\def\NSoplus{\mbox{$\bigcirc\hspace{-0.122in}+\hspace{0.055in}$}}

\def\NSoone{\mbox{$\bigcirc\hspace{-0.095in}{\scriptstyle{1}}\hspace{0.055in}$}}
\def\NStinyoplus{\tiny \mbox{$\bigcirc\hspace{-0.085in}+\hspace{0.012in}$}}
\def\NStinyotimes{\tiny \mbox{$\bigcirc\hspace{-0.085in}\times\hspace{0.012in}$}}
\def\tinyoplus{\tiny \mbox{$\hspace{0.015in}\bigcirc\hspace{-0.085in}+\hspace{0.012in}$}}
\def\tinyotimes{\tiny \mbox{$\hspace{0.015in}\bigcirc\hspace{-0.085in}\times\hspace{0.012in}$}}
\makeindex

\title{A Graphic Generalization of Arithmetic}
\author{Bilal Khan
\thanks{Department of Mathematics and Computer Science, City University of New York John Jay College of Criminal Justice, 899 Tenth Avenue, New York, NY 10019
({\tt grouptheory@hotmail.com}).}
\and
Kiran R. Bhutani
\thanks{Department of Mathematics, The Catholic University of America, Washington DC 20064
({\tt bhutani@cua.edu}).}
\and
Delaram Kahrobaei
\thanks{Department of Mathematics, City University of New York, Graduate Center, New York, NY 10016
({\tt dkahrobaei@gc.cuny.edu}).}}

\begin{document}
\maketitle
%-----------------------------------------------------------------------
%-----------------------------------------------------------------------
\begin{abstract}
In this paper, we extend the classical arithmetic defined over the set
of natural numbers ${\mathbb N}$, to the set of all finite directed
% MOD2
connected multigraphs having a pair of distinguished vertices.
Specifically, we introduce a model $\mathcal F$ on the set of such
graphs, and provide an interpretation of the language of arithmetic
${\mathcal L} = \{0,1,\leqslant,+,\times\}$ inside $\mathcal F$. The
resulting model exhibits the property that the standard model on
${\mathbb N}$ embeds in $\mathcal F$ as a submodel, with the directed
path of length $n$ playing the role of the standard integer $n$. We
will compare the theory of the larger structure $\mathcal F$ with
classical arithmetic statements that hold in ${\mathbb N}$.  For
example, we explore the extent to which $\mathcal F$ enjoys properties
like the associativity and commutativity of $+$ and $\times$,
distributivity, cancellation and order laws, and decomposition into
irreducibles.
\end{abstract}

\begin{keywords} 
arithmetic, graphs.
\end{keywords}

\begin{AMS}
05C99, 11U10
\end{AMS}

\pagestyle{myheadings}
\thispagestyle{plain}

\markboth{BILAL KHAN, KIRAN R. BHUTANI AND DELARAM KAHROBAEI}{A GRAPHIC GENERALIZATION OF ARITHMETIC}

%-----------------------------------------------------------------------
%-----------------------------------------------------------------------
\section{Introduction}

The {\em language of arithmetic} ${\mathcal L}$ consists of two
$0$-ary relations ${\mathbf 0}$ and ${\mathbf 1}$, one binary relation
${\mathbf \leqslant}$, and two ternary relations ${\mathbf +}$ and
${\mathbf \times}$.  In this paper, we generalize classical arithmetic
defined over the natural numbers ${\mathbb N} = \{0, 1, 2, \ldots \}$,
% MOD2
to the set $F$ consisting of all {\em flow graphs}: finite directed
connected multigraphs in which a pair of distinguished vertices designated
as the {\em source} and {\em target} vertex.  We give natural
interpretation for ${\mathcal L}$ on the set $F$. To avoid confusion
with the standard model of arithmetic, the corresponding operations in
$F$ are denoted with a circumscribed circle.  The new model ${\mathcal
F} = \langle F, \ozero, \oone, \weaklesseq, \oplus, \otimes \rangle$
is a natural extension of the standard model ${\mathcal N} = \langle
{\mathbb N}, 0, 1, \leqslant +, \times \rangle$.  Specifically, we
exhibit an embedding $i: {\mathcal N} \stackrel{i}{\hookrightarrow}
{\mathcal F}$ satisfying:

\begin{eqnarray*}
i(0) & = & \ozero,\\
i(1) & = & \oone,\\
\forall x,y \in {\mathbb N}, \;\; x \leqslant y & \Rightarrow & i(x) \weaklesseq i(y),\\
\forall x,y \in {\mathbb N}, \;\; i(x+y) & = & i(x) \oplus i(y),\\
\forall x,y \in {\mathbb N}, \;\; i(x\times y) & = & i(x) \otimes i(y).\\
\end{eqnarray*}

{\bf Objective:} {\em Compare the theory $Th({\mathcal F}) = \{ \phi
\;|\; {\mathcal F} \models \phi \}$ with {\em true arithmetic} $TA =
\{ \phi \;|\; {\mathcal N} \models \phi \}$}\footnote{Following
standard model theory, here $\phi$ is a first-order sentence in the
language ${\mathcal L}$.}.

\vspace{0.1in} 

There have been other attempts to define algebraic and metric
structures on the set of all graphs.  The classical operations on
graphs \cite{west} (including extensive literature on graph products
\cite{imrichbook}) have yielded deep results and a profound
mathematical theory. However, to date, these operations have not
provided an interpretation of the language of arithmetic on graphs.
This paper presents results and open questions in this direction.  In
\cite{BK3,BK2,BK1}, the authors used graph embeddings to define a
metric on the set of all simple connected graphs of a given order.
This work differs from those investigations in that it considers an
infinite collection of graphs in order to extend the standard model of
arithmetic, and in doing so does not seek to establish a metric
structure.

% \nocite{GP}
% One notable instance along these lines is the work of Kontsevich on
% graph homology \cite{K2,K1}, which was subsequently used to compute
% the homology of a certain infinite dimensional Lie algebras and to
% parameterize invariants of certain odd dimensional manifolds.  
%-----------------------------------------------------------------------

\begin{definition}[Flow graph]
\label{Flow graph}
\index{flow graph}
\index{flow graph!finite}
\index{flow graph!infinitesimal}
\index{finite!flow graph}
\index{infinitesimal!flow graph}
\index{Graph!flow}
\index{flow graph!isomorphism}
\index{isomorphism!flow graph}
% MOD2
A {\em flow graph} $A$ is a triple $(G_A, s_A, t_A)$, where $G_A$ is a
finite directed connected multigraph and $s_A, t_A \in V[G_A]$ are
called the {\em source} and the {\em target} vertex of $A$,
respectively.  The set of all flow graphs is denoted $F$.  The unique
flow graph for which $|V[G_A]| = 1$ and $|E[G_A]| = 0$ is called {\em
the trivial flow graph}; all other flow graphs are considered {\em
non-trivial}.  If $s_A=t_A$ and $A$ is non-trivial then $A$ is called
an {\em infinitesimal} flow graph.  The set of all infinitesimal flow
graphs is denoted $I$. 
% REV IM
Given two flow graphs $A=(G_A, s_A, t_A)$ and $B=(G_B, s_B, t_B)$, a
map $\phi: A \rightarrow B$ is called an injective morphism of flow
graphs if (as a graph embedding) $\phi$ maps $G_A$ injectively into
$G_B$ and additionally satisfies $\phi(s_A) = s_B$, $\phi(t_A) = t_B$.
Flow graphs $A$ and $B$ are considered isomorphic if there is an
injective morphism $\phi: A \rightarrow B$ for which $Im(\phi) = B$.
\end{definition}

\begin{definition}[Graphical natural number]
\label{Graphical natural number} 
\index{Graphical natural number}
\index{Natural number!graphical} 
We represent the natural number $n$ as a directed chain of length $n$,
having $n+1$ vertices. More formally, let $P_n$ be a directed chain of
length $n$ (having $n+1$ vertices) where each vertex has in-degree
$\leqslant 1$ and out-degree $\leqslant 1$.  Denote by $s_n$, the
unique vertex in $P_n$ having in-degree $0$, and let $t_n$ be the
unique vertex in $P_n$ having out-degree $0$.  The flow graph
$F_n=(P_n,s_n,t_n)$ is referred {\em the graphic natural number
$n$}. Define the map $i:{\mathcal N} \rightarrow {\mathcal F}$ as
\begin{eqnarray*}
  i: n & \mapsto & F_n.
\end{eqnarray*}
\end{definition}

%-----------------------------------------------------------------------
\subsection{Addition}

\begin{figure}[htb]
\centering{\mbox{\psfig{figure=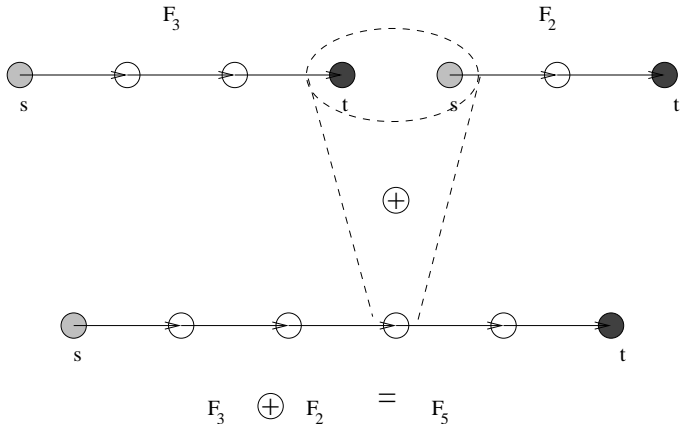}}}
\caption{Interpreting addition of natural numbers inside
${\mathcal F}$.} \label{plus-standard-picture}
\end{figure}

In Definition~\ref{Flow graph}, we represented the natural number
$n$ by the flow graph $F_n$.  It follows that we interpret the
addition of two numbers $n_1$ and $n_2$ inside ${\mathcal F}$ as
``concatenating'' $F_{n_1}$ with $F_{n_2}$.  Consider, for
example, the addition of $3$ and $2$ depicted in
Figure~\ref{plus-standard-picture}.

\begin{figure}[htb]
\centering{\mbox{\psfig{figure=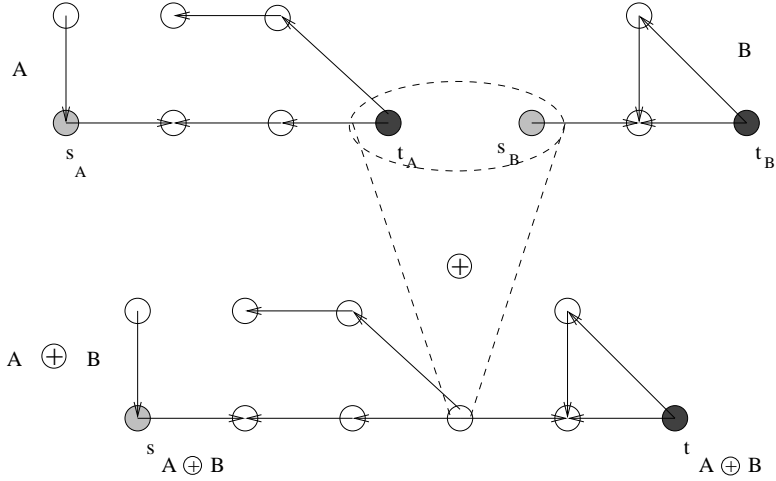}}}
\caption{General addition of flow graphs.}
\label{plus-general-picture}
\end{figure}

To extend this definition of $\oplus$ to all of $F$, we define general
addition of flow graphs as follows: Given two flow graphs $A$ and $B$,
define $A\oplus B$ to be the flow graph obtained by identifying $t_A$
with $s_B$ and defining $s_{A\tinyoplus B} = s_A$ and $t_{A\tinyoplus B} =
t_B$.  An example of such an addition is shown in
Figure~\ref{plus-general-picture}.

To make this formal we define the following operation on connected
directed multigraphs: Given directed graphs $G_1$ and $G_2$, and
vertices $u_1 \in V[G_1]$, $u_2 \in V[G_2]$, we define
\begin{eqnarray*}
  G_1 \oplus_{u_1 \approx u_2} G_2 & \definemath & (G_1 \sqcup G_2) /
(u_1 \approx u_2)
\end{eqnarray*}
to be the graph obtained by taking disjoint copies of $G_1$ and $G_2$
and identifying vertex $u_1$ in $G_1$ with vertex $u_2$ in $G_2$.
Note the obvious and natural injective graph homomorphisms
\begin{eqnarray}
\begin{array}{rcl}
  \sigma^{\NStinyoplus}_{u_1 \approx u_2} : G_1 & \hookrightarrow & G_1
\oplus_{u_1 \approx u_2} G_2 \\
  \tau^{\NStinyoplus}_{u_1 \approx u_2} : G_2 & \hookrightarrow & G_1
\oplus_{u_1 \approx u_2} G_2.
\end{array}
\label{plus-injections}
\end{eqnarray}

%%%%%%%%%%%%%%%%%%
\begin{definition}
\label{oplus-definition}
\index{flow graph!addition}
Given two flow graphs $A=(G_A, s_A, t_A)$ and $B=(G_B, s_B, t_B)$, we
define
\begin{eqnarray*}
  A \oplus B & \definemathit & ( G_A \oplus_{t_A\approx s_B} G_B, s_A,
t_B).
\end{eqnarray*}
\end{definition}

\begin{remark}
\label{plus-remark}
% MOD2
Note that if $A$ is a flow graph with $p_A$ vertices and $q_A$ edges,
and $B$ is a flow graph with $p_B$ vertices and $q_B$ edges, then $A
\oplus B$ is a flow graph having $p_A + p_B - 1$ vertices and $q_A +
q_B$ edges.
\end{remark}

%%%%%%%%%%%%%%%%%%
% MOD2
The next lemma follows immediately from Definitions
\ref{Graphical natural number} and \ref{oplus-definition}.

\begin{lemma}
\label{plus-embeds}
Let $m, n$ be natural numbers.  Then $i(n+m) = i(n) \oplus i(m)$.
\end{lemma}

\begin{lemma}
\label{plus-identity} \index{$\oplus$!identity} \index{$\ozero$}
$\ozero \definemathit F_0$ is the unique two-sided identity with
respect to $\oplus$.  That is, for all flow graphs $A,G \in {\mathcal
F}$,
\begin{eqnarray*}
\begin{array}{rcccl}
     A \oplus G = A & \Leftrightarrow & G = \ozero & \Leftrightarrow & G \oplus A = A.
\end{array}
\end{eqnarray*}
\end{lemma}
\begin{proof}
% MOD2
If $G=F_0$ then $A \oplus G = G \oplus A = A$.  For the reverse, we
appeal to Remark~\ref{plus-remark}, noting that $A \oplus G = A$
implies $p_A + p_G - 1 = p_A$ and $q_A + q_G = q_A$.  Hence $p_G=1$
and $q_G=0$, so $G=F_0$.  An analogous argument shows that $G \oplus A
= A$ implies $G=F_0$.
\end{proof}

\begin{observation}
We note that the sum of two infinitesimals is again an infinitesimal.
On the other hand, if at least one summand is a non-trivial
non-infinitesimal flow graph then the summation evaluates to a
non-trivial non-infinitesimal flow graph.
\end{observation}

\begin{definition}[Scalar multiplication of flow graphs]
Given a flow graph $A$, and a positive natural number $k$ in ${\mathbb
N}$, we define left-multiplication inductively as follows:
\begin{eqnarray*}
     1A & = & A \\
     kA & = & (k-1)A \oplus A.
\end{eqnarray*}
Right-multiplication is defined analogously.  However, as we will see,
$\oplus$ is associative, and so the two notions coincide.  We shall
subsequently consider only left-multiplication by integer scalars.
\label{Scalar multiplication of flow graphs}
\index{Scalar multiplication of flow graphs}
\index{flow graph!scalar multiplication}
\end{definition}

%-----------------------------------------------------------------------

\subsection{Multiplication}

In the previous section, we presented an interpretation of
addition in ${\mathcal F}$ that is a natural extension of addition
on the natural numbers.  In this section, we give an
interpretation of multiplication in ${\mathcal F}$.  In doing
this, we must respect the fact that for each pair of natural
numbers $n_1, n_2$, the following identity holds in ${\mathcal
N}$:

\begin{eqnarray*}
I_{n_1, n_2}: \underbrace{n_2 + n_2 + \cdots + n_2}_{n_1\text{ times}} =
n_1 n_2 = \underbrace{n_1 + n_1 + \cdots + n_1}_{n_2\text{ times}}.
\end{eqnarray*}
So, in particular, the definition of $\otimes$ in ${\mathcal F}$
must satisfy
\begin{eqnarray}
\label{multiplication-addition-relation}
n_1 F_{n_2} = F_{n_1} \otimes F_{n_2} = n_2 F_{n_1}.
\end{eqnarray}

\begin{figure}[htb]
\centering{\mbox{\psfig{figure=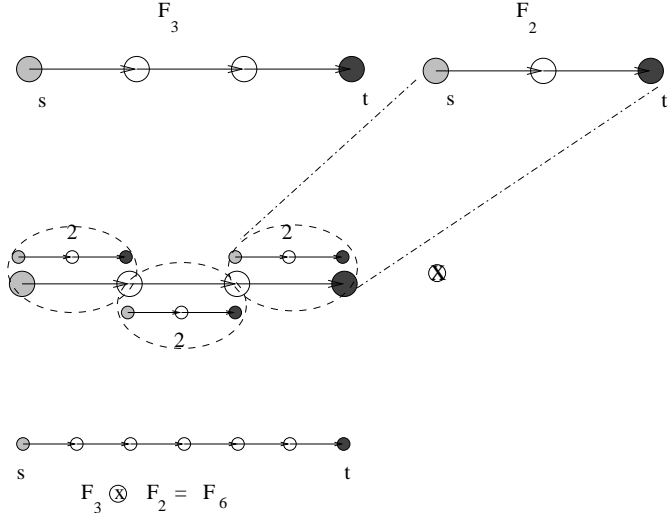}}}
\caption{Standard multiplication of natural numbers in $\mathcal{F}$
(represented as flow graphs).}
\label{times-standard-picture}
\end{figure}

Given that we represent the natural number $n$ by the flow graph
$F_n$, the product of two graphical numbers $F_{n_1}$ and $F_{n_2}$
can be made to satisfy
relation~(\ref{multiplication-addition-relation}) if we take
multiplication to be the act of replacing each edge of $F_{n_1}$ with
a copy of $F_{n_2}$.  Consider the multiplication of graphical natural
numbers $F_3$ and $F_2$, as depicted in
Figure~\ref{times-standard-picture}.

\begin{figure}[htb]
\centering{\mbox{\psfig{figure=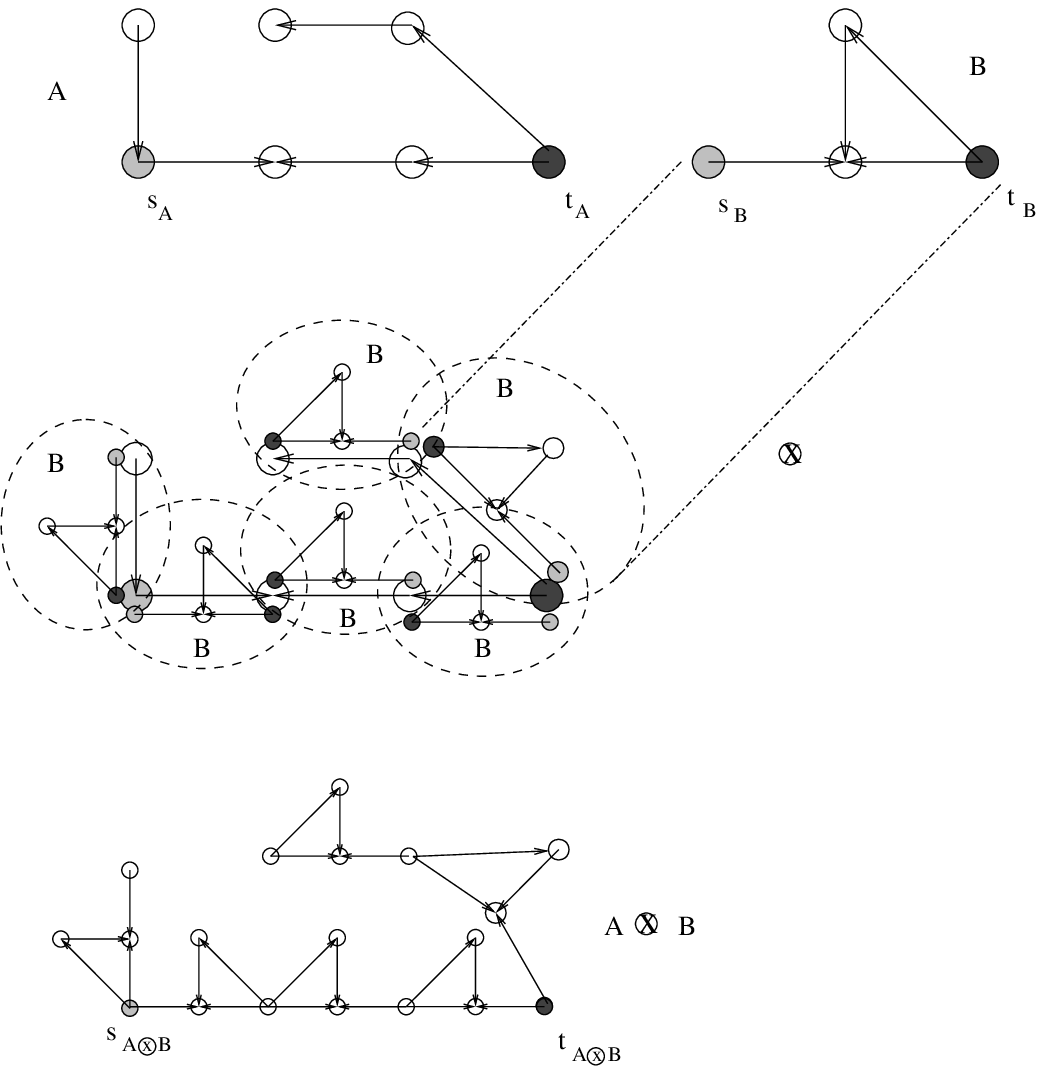}}}
\caption{General multiplication of flow graphs.}
\label{times-general-picture}
\end{figure}

To extend this definition of $\otimes$ to all of $F$, we define
general multiplication of flow graphs as follows: Given two flow
graphs $A$ and $B$, define $A\otimes B$ to be the flow graph obtained
by replacing every edge $e$ (from $E[G_A]$) with a copy of $B$ as
follows: For each edge $e=(u,v)$ in $A$, we remove $e$ and replace it
with a graph $B_e$ isomorphic to $B$, by identifying $u$ with
$s_{B_e}$, and $v$ with $t_{B_e}$.  An example of such a
multiplication is shown in Figure~\ref{times-general-picture}.

To make this formal we define the following operation on connected
directed multigraphs: Given directed graphs $G_1$ and $G_2$, an edge
$e=(u_1,v_1) \in E[G_1]$ and vertices $u_2, v_2 \in V[G_2]$, we define
\begin{eqnarray*}
  G_1 \otimes_{e \approx (u_2,v_2)} G_2 & \definemath & [(G_1 \backslash
e) \sqcup G_2] / (u_1 \approx u_2, v_1 \approx v_2)
\end{eqnarray*}
to be the graph obtained by removing $e$ from $G_1$ and attaching
a copy of $G_2$ to the resulting graph by gluing $u_1$ with $u_2$
and $v_1$ with $v_2$.
Note the obvious and natural injective maps
\begin{eqnarray}
\begin{array}{rcl}
  \sigma^{\NStinyotimes}_{e \approx (u_2,v_2)} : G_1\backslash e &
\hookrightarrow & G_1 \otimes_{e \approx (u_2,v_2)} G_2 \\
  \tau^{\NStinyotimes}_{e \approx (u_2,v_2)} : G_2 & \hookrightarrow & G_1
\otimes_{e \approx (u_2,v_2)} G_2.
\end{array}
\label{times-injections}
\end{eqnarray}

%%%%%%%%%%%%%%%%%%
\begin{definition}
\label{otimes-definition}
\index{flow graph!multiplication} Given flow graphs $A=(G_A, s_A,
t_A)$ and $B=(G_B, s_B, t_B)$.  Define the directed graph $K_0 = G_A$,
$E_0 = E[G_A]$, and let $\bar{\sigma}_0, \sigma_0:K_0 \rightarrow K_0$
be the identity isomorphisms.  Fix any enumeration $\eta$ of the edges
$E[G_A]$, say $\eta = e_1, e_2, \ldots, e_m$.  Inductively, for $i=1,
2, \ldots, m$ we define
\begin{eqnarray*}
  K_i & = & K_{i-1} \otimes_{\bar{\sigma}_{i-1}(e_i) \approx (s_B, t_B)}
G_B \\
  \sigma_i & = & \sigma^{\NStinyotimes}_{\bar{\sigma}_{i-1}(e_i) \approx
(s_B, t_B)} : K_{i-1} \backslash e_i \rightarrow K_i \\
  E_i & = & E_{i-1} \backslash \{e_i\} \\
  \bar{\sigma}_{i} & = &
(\sigma_{i}){}_{|E_i}(\sigma_{i-1}){}_{|E_i}\cdots
(\sigma_{1}){}_{|E_i}(\sigma_{0}){}_{|E_i}.
\end{eqnarray*}
Informally, $K_i$ is the directed graph obtained after edges $e_1,
\ldots, e_i$ have been deleted from $G_A$ and replaced by copies of
$G_B$.  Finally, we put
\begin{eqnarray*}
  A \otimes_\eta B & \definemath & ( K_m, \bar{\sigma}_{m}(s_A),
\bar{\sigma}_{m}(t_A)).
\end{eqnarray*}
\end{definition}

The reader may verify that the operation $\otimes_\eta$ is
well-defined, and that in particular, it is independent of the chosen
enumeration $\eta$ of the edges $E[G_A]$.

\begin{remark}
\label{times-remark}
% MOD3
Let $A$ be a flow graph with $p_A$ vertices and $q_A$ edges, and $B$
be a flow graph having $p_B$ vertices and $q_B$ edges.  Then $A
\otimes B$ has $q_A q_B$ edges.  If $B$ is either trivial or
infinitesimal then $A \otimes B$ has $1 + q_A(p_B - 1)$ vertices.  
% REV IM
If $B$ is non-trivial and not infinitesimal then $A \otimes B$ has
$p_A + q_A(p_B - 2)$ vertices.
% MOD3
% If $A$ and $B$ are both non-trivial and not infinitesimal then $A
% \otimes B$ has $p_A + q_A(p_B - 2)$ vertices.
% 
% If $A$ is non-trivial and not infinitesimal, and $B$ is either trivial
% or infinitesimal then $A \otimes B$ has $1 + q_A(p_B - 1)$ vertices.
% 
% If $A$ is either trivial or infinitesimal, and $B$ is non-trivial and
% not infinitesimal then $A \otimes B$ has $p_B + q_B(p_B - 2)$
% vertices.
% 
% If both $A$ and $B$ are either trivial or infinitesimal, then $A
% \otimes B$ has $1 + q_A(p_B - 1)$ vertices.
\end{remark}

%%%%%%%%%%%%%%%%%%
% MOD2
The next lemma follows immediately from Definitions 
\ref{Graphical natural number} and \ref{otimes-definition}.

\begin{lemma}
\label{times-embeds}
Let $m, n$ be natural numbers.  Then $i(n\times m) = i(n) \otimes
i(m)$.
\end{lemma}

\begin{lemma}
\label{times-identity} 
\index{$\otimes$!identity} 
\index{$\oone$}
% MOD2
Given flow graphs $G$ and $H$, and a non-trivial, non-infinitesimal
flow graph $A$:
\begin{eqnarray*}
&
\begin{array}{rcccl}
A \otimes G = A   & \Leftrightarrow & G = \oone & \Leftrightarrow & G \otimes A = A,\\
G \otimes H = \ozero & \Leftrightarrow & H = \ozero & \text{ or  } & G = \ozero.
\end{array}&
\end{eqnarray*}
\end{lemma}
\begin{proof}
% MOD3
If $G=F_1$ then $A \otimes G = G \otimes A = A$.  For the reverse, we
appeal to Remark~\ref{times-remark}, noting that $A \otimes G = A$
implies $p_A + q_A(p_G - 2) = p_A$ and $q_A q_G = q_A$.  Hence
$p_G=2$ and $q_G=1$, so $G=F_1$.  An analogous argument shows that $G
\otimes A = A$ implies $G=F_1$.  

If $G=F_0$ then $G \otimes H = H \otimes G = F_0$.  For the reverse, we
appeal to Remark~\ref{times-remark}, noting that $G \otimes H =
\ozero$ implies $q_G q_H = 0$, so either $q_G=0$ or $q_H=0$.  It
follows that either $H = \ozero$ or $G = \ozero$.
\end{proof}

% MOD3
Note that if we remove the hypothesis that $A$ is non-infinitesimal in
Lemma \ref{times-identity}, then $A \otimes G = A$ and $G \otimes A =
A$ do not necessarily imply $G = \oone$.  The simplest counterexample
is seen by taking $G=A$ to be the flow graph consisting of one vertex
and one loop edge.  We denote this graph as $C_1$.  Indeed, it is easy
to see that this is the only counterexample.  Suppose $A \otimes G =
A$, for some infinitesimal $A$.  By Remark \ref{times-remark}, this
implies $1 + q_A(p_G - 1) = p_A$ and $q_A q_G = q_A$.  Hence $q_G=1$
and $p_G=1 + (p_A - 1)/q_A$.  If $q_G=1$ then either $p_G = 2$ or $1$.
If $p_G=2$, then $G=F_1$, and this is not a counterexample.  If
$p_G=1$, then $G$ is a graph consisting of one vertex and one loop
edge.  Moreover, since $q_A$ is finite, $1 + q_A(p_G - 1) = p_A$
implies that $p_A=1$.  It follows that $A = q_A C_1$.

\begin{observation}
Suppose $G$ and $H$ are non-trivial flow graphs of which {\em at least
one}, say $H$, is infinitesimal.  Then $s_{G \tinyotimes H}=t_{G
\tinyotimes H}$ in $G \otimes H$ and $s_{H \tinyotimes G}=t_{H
\tinyotimes G}$ in $H \otimes G$.  Hence $G \otimes H$ and $H \otimes
G$ are both infinitesimal.

On the other hand, suppose $G$ and $H$ are non-trivial flow graphs that are
both non-infinitesimal.  Then $s_{G \tinyotimes H}\neq t_{G
\tinyotimes H}$ in $G \otimes H$ and $s_{H \tinyotimes G}\neq t_{H
\tinyotimes G}$ in $H \otimes G$.  Hence $G \otimes H$ and $H \otimes
G$ are both non-infinitesimal.

It follows that if $G$ and $H$ are non-trivial flow graphs, then $G
\otimes H$ is infinitesimal if and only if at least one of the two
factors is infinitesimal.  The reader may wish to compare this with
the second assertion of Lemma \ref{times-identity}.
\end{observation}

\begin{definition}[Scalar exponentiation of flow graphs]
\label{Scalar exponentiation of flow graphs}
\index{Scalar exponentiation of flow graphs}
\index{flow graph!scalar exponentiation}
Given a flow graph $A$, and a positive natural number $k$ in ${\mathbb
N}$, we define right-exponentiation inductively as follows:
\begin{eqnarray*}
  A^1 & = & A\\
  A^k & = & A^{k-1} \otimes A.
\end{eqnarray*}
Left-exponentiation is defined analogously.  However, as we will see
shortly, $\otimes$ is associative, and so the two notions coincide.
We shall subsequently consider only right-exponentiation by integer
scalars.
\end{definition}

%-----------------------------------------------------------------------
\subsection{Order}

Given our representation of the natural number $n$ by the flow graph
$F_n$ in Definition~\ref{Graphical natural number}, comparing the
order of two numbers $n_1$ and $n_2$ amounts to simply comparing the
lengths of the corresponding chain graphs $F_{n_1}$ and $F_{n_2}$.  To
generalize this to all of ${\mathcal F}$, however, we cannot refer to
``length''.  In what follows, we present two possible interpretations
of $\leqslant$ in ${\mathcal F}$.  To avoid confusion, we refer to
these interpretations as $\weaklesseq$ and $\stronglesseq$.
%--------------------------------------

\subsubsection{Weak Order $\weaklesseq$}

Suppose we are given two flow graphs $A$ and $B$.  Informally, we
say that $A \weaklesseq B$ iff there is a way to partition $A$ into
edge-disjoint neighborhoods of the source/target of vertices of $A$ in
such a way that these neighborhoods can be mapped into disjoint
neighborhoods of the source/target vertices of $B$.
To make this more precise we define the following operation on
connected directed multigraphs.

\begin{definition}[$(s,t)$-splitting]
\label{(s,t)-splitting}
\index{$(s,t)$-splitting}
\index{splitting!$(s,t)$}
Given a connected directed multigraph $G=(V,E)$ and two vertices $s$
and $t$ in $V$, an $(s,t)$-{\em splitting} of $G$ is a pair of graphs
$(H_1, H_2)$ with the following properties:
\begin{itemize}
\item $H_1$ and $H_2$ are connected subgraphs of $G$.
\item $s$ is in $V[H_1]$ and $t$ is in $V[H_2]$.
\item $\{E[H_1], E[H_2]\}$ is a partition of $E$.  While this implies
$V[H_1] \cup V[H_2] = V[G]$, we remark that $V[H_1] \cap V[H_2]$ need
not be empty.
\end{itemize}
\end{definition}

We can now give a precise definition of the weak ordering.
\begin{definition}[Weak order]
\label{weak order}
\index{weak order}
\index{order!weak}
Given two flow graphs $A=(G_A, s_A, t_A)$ and $B=(G_B, s_B, t_B)$, we say
that $A \weaklesseq B$ if there is an $(s_A,t_A)$-splitting $(H_1,
H_2)$ of $G_A$ and graph embeddings $\phi_1:H_1 \rightarrow G_B$,
$\phi_2:H_2 \rightarrow G_B$ such that $\phi_1(s_A) = s_B$ and
$\phi_2(t_A) = t_B$ and $\phi_1(E[H_1]) \cap \phi_2(E[H_2]) =
\emptyset$.
\end{definition}

\begin{figure}[htb]
\centering{\mbox{\psfig{figure=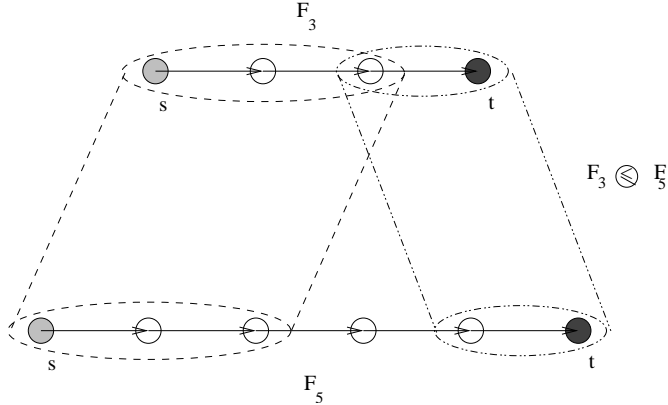}}}
\caption{Standard weak ordering of natural numbers (represented as
flow graphs).}
\label{order-standard-picture}
\end{figure}

Consider the comparison of $F_3$ and $F_5$ in
Figure~\ref{order-standard-picture} which illustrates the assertion that
$F_3
\weaklesseq F_5$.

The proof of the following lemma is immediate.

\begin{lemma}
\label{order-embeds}
Let $m, n$ be natural numbers.  Then $n \leqslant m \Leftrightarrow
i(n) \weaklesseq i(m)$.
\end{lemma}

\begin{figure}[htb]
\centering{\mbox{\psfig{figure=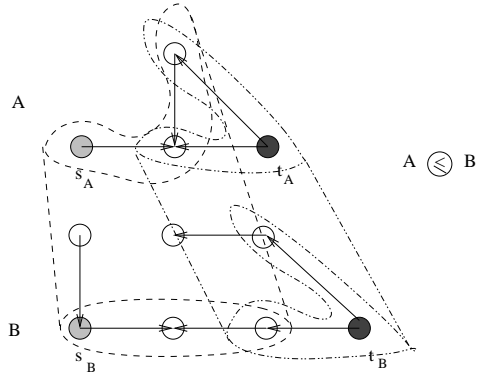}}}
\caption{General weak ordering of flow graphs.}
\label{order-general-picture}
\end{figure}

Figure~\ref{order-general-picture} illustrates a more general
example in which weak order is used to compare two elements of
${\mathcal F}$ which are {\em not} graphical natural numbers.

The next Proposition follows immediately from Lemmas~\ref{plus-embeds},
\ref{plus-identity}, \ref{times-embeds},
\ref{times-identity}, and \ref{order-embeds}.

\begin{proposition}
\label{submodel} 
Under the embedding $i:n \mapsto F_n$, the standard model ${\mathcal
N} = \langle {\mathbb N}, 0, 1, \leqslant, +, \times \rangle$ is a submodel of
${\mathcal F} = \langle F, \ozero, \oone, \weaklesseq, \oplus, \otimes
\rangle$, where $\ozero = F_0$, $\oone = F_1$, and the relations
$\oplus$, $\otimes$ and $\weaklesseq$ reinterpret $+, \times$ and $\leqslant$
inside ${\mathcal F}$.
\end{proposition}

%--------------------------------------
\subsubsection{Strong Order $\stronglesseq$}

We now give an alternate, strengthened ordering on ${\mathcal F}$.
Given two flow graphs $A$ and $B$, informally, we say that $A
\stronglesseq B$ iff a copy of $G_A$ appears as a neighborhood of
both $s_B$ and $t_B$ in $G_B$.  The next definition makes this
statement precise.

\begin{definition}[Strong order]
\label{strong order}
\index{strong order}
\index{order!strong}
Given two flow graphs $A=(G_A, s_A, t_A)$ and $B=(G_B, s_B, t_B)$, we
say $A \stronglesseq B$ iff there are graph embeddings
$\phi_s:G_A \rightarrow G_B$ and $\phi_t:G_A \rightarrow G_B$ which
satisfy $\phi_s(s_A) = s_B$ and $\phi_t(t_A) = t_B$.
\end{definition}

\begin{figure}[htb]
\centering{\mbox{\psfig{figure=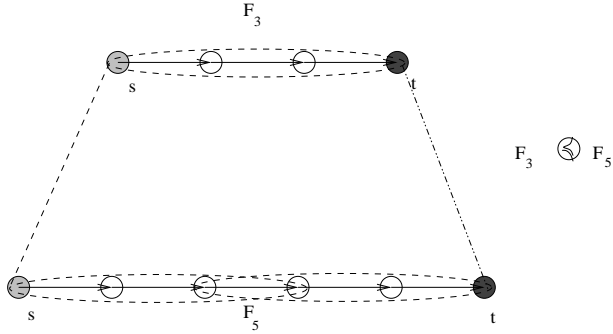}}}
\caption{Standard strong ordering of natural numbers (represented as
flow graphs).}
\label{order-standard-picture-strong}
\end{figure}

Consider the comparison of $F_3$ and $F_5$ depicted in
Figure~\ref{order-standard-picture-strong}; clearly $F_3 \stronglesseq
F_5$.

The proof of the following lemma is immediate.

\begin{lemma}
\label{order-embeds-strong}
Let $m, n$ be natural numbers.  Then $n \leqslant m \Leftrightarrow
i(n) \stronglesseq i(m)$.
\end{lemma}

\begin{figure}[htb]
\centering{\mbox{\psfig{figure=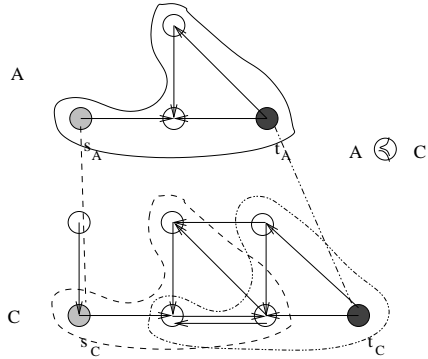}}}
\caption{General strong ordering of flow graphs.}
\label{order-general-picture-strong}
\end{figure}

Figure~\ref{order-general-picture-strong} illustrates a more
general example in which strong order is used to compare two
elements of ${\mathcal F}$ which are {\em not} graphical natural
numbers.

The next Proposition follows immediately from Lemmas~\ref{plus-embeds},
\ref{plus-identity}, \ref{times-embeds},
\ref{times-identity}, and \ref{order-embeds-strong}.
\begin{proposition}
\label{submodel-strong} Under the embedding $i:n \mapsto F_n$, the
standard model ${\mathcal N} = \langle {\mathbb N}, 0, 1, \leqslant, +,
\times \rangle$ is a submodel of ${\mathcal F} = \langle F,
\ozero, \oone, \stronglesseq, \oplus, \otimes \rangle$, where
$\ozero = F_0$, $\oone = F_1$, and the relations $\oplus$,
$\otimes$ and $\stronglesseq$ reinterpret $+, \times$ and $\leqslant$ inside
${\mathcal F}$.
\end{proposition}

The next proposition and example show that ordering by $\stronglesseq$ is indeed
strictly stronger than ordering by $\weaklesseq$.
\begin{proposition}
\label{strong-order-is-stronger}
Given flow graphs $A=(G_A, s_A, t_A)$ and $B=(G_B, s_B, t_B)$
\begin{eqnarray*}
A \stronglesseq B & \Rightarrow & A \weaklesseq B.
\end{eqnarray*}
\end{proposition}
\begin{proof}
Since $A \stronglesseq B$, there are graph embeddings $\phi_s:G_A
\rightarrow G_B$ and $\phi_t:G_A \rightarrow G_B$ which satisfy
$\phi_s(s_A) = s_B$ and $\phi_t(t_A) = t_B$.
Let $E_1 = E[G_A]$, $V_1 = V[G_A]$; take $E_2 = \emptyset$, $V_2 =
\{t_A\}$.  Put $H_1 = (V_1, E_1)$ and $H_2 = (V_2, E_2)$.  Then $(H_1,
H_2)$ is an $(s_A,t_A)$-splitting of $G_A$.  We take graph embeddings
$\phi_1 = \phi_s|_{H_1}:H_1 \rightarrow G_B$, and $\phi_2 = \phi_t|_{H_2}:H_2
\rightarrow G_B$.  Then $\phi_1(s_A) = s_B$ and $\phi_2(t_A) = t_B$
and $\phi_1(E[H_1]) \cap \phi_2(E[H_2]) = \emptyset$.  Thus, $A
\weaklesseq B$.
\end{proof}

The converse of Proposition \ref{strong-order-is-stronger} is false,
as the following example indicates.  

\begin{example}
Take $A$ and $B$ to be the flow graphs depicted on page
\pageref{order-general-picture}, where
Figure~\ref{order-general-picture} illustrates that $A \weaklesseq B$.
Note that $G_A$ contains a vertex of degree $3$, while
$G_B$ does not, hence no neighborhood of $s_B$ or $t_B$ can be
isomorphic to $G_A$.  Thus $A \notstronglesseq B$.
\end{example}
%-----------------------------------------------------------------------
%-----------------------------------------------------------------------

\section{Results}

% MOD2: STILL NEED TO DO THIS!!!
We begin by considering properties of $\oplus$ in
Section~\ref{addition-section}.  We show that $\oplus$ is an
associative, non-commutative operation, and provide a natural
criterion for a flow graph to be irreducible as a proper sum.  We
prove that every flow graph is canonically decomposable as a sum of
irreducibles. Using this canonical decomposition, we deduce left and
right cancellation laws for $\oplus$, and show that if two flow graphs
$A$ and $B$ commute with respect to $\oplus$ then they are necessarily
scalar multiples of some flow graph $C$.  Then, in
Section~\ref{multiplication-section} we show that $\otimes$ is an
associative, non-commutative operation and that it right-distributes
over $\oplus$ (but does not left-distribute). We define left and right
divisibility of flow graphs, and use this to introduce the notion of a
prime flow graph, and show that the concept of left-prime and
right-prime coincide. We describe the canonical $\oplus$ decomposition
of flow graph products in terms of the $\oplus$ decompositions of each
of the $\otimes$ factors.  Finally, in Section~\ref{order-section}, we
explore the relationship between strong order (denoted by
$\stronglesseq$) and weak order (denoted by $\weaklesseq$), describing
the interaction between these orders and the operations of $\oplus$
and $\otimes$.  We show that while the two orders coincide on the
graphical natural numbers, neither order is anti-symmetric on all of
${\mathcal F}$, and only $\stronglesseq$ is transitive.  On the other
hand, many of the laws that govern the relationship between
$\leqslant$, $+$ and $\times$ in ${\mathcal N}$ continue to hold for
$\weaklesseq$, $\oplus$ and $\otimes$ in ${\mathcal F}$, but these
laws are violated under the ordering $\stronglesseq$.

%------------------------------------------------------
\subsection{Additive Properties}
\label{addition-section}

In this section we present some properties of $\oplus$.

\begin{lemma}[Associativity of $\oplus$]
\label{oplus associative}
\index{$\oplus$!associativity}
The operation $\oplus$ is associative.
\end{lemma}
\begin{proof}
Given flow graphs $A, B, C$,
\begin{eqnarray*}
(A \oplus B) \oplus C & = & ( G_{A} \oplus_{t_A\approx s_B} G_{B}, s_A,
t_B) \oplus C\\
                                                & = & ((G_{A}
\oplus_{t_A\approx s_B} G_{B}) \oplus_{t_B\approx s_C} G_{C}, s_A,
t_C)\\
                                                & = & (G_{A}
\oplus_{t_A\approx s_B} (G_{B} \oplus_{t_B\approx s_C} G_{C}), s_A,
t_C)\\
                                                & = & A \oplus ( G_{B}
\oplus_{t_B\approx s_C} G_{C}, s_B, t_C)\\
                                                & = & A \oplus (B \oplus
C).
\end{eqnarray*}
\end{proof}

\begin{example}
\label{non-comm-plus-example}
Let $A$ be the flow graph consisting of a directed cycle of length $3$
and let source and target vertices be any two vertices on this cycle.
Then it is easy to check that $A \oplus F_2$ is not equal to $F_2
\oplus A$, that is to say, there is no flow graph isomorphism between
$A\oplus F_2$ and $F_2\oplus A$ (see Figure
\ref{non-comm-plus-figure}).  
\begin{figure}[htb]
\centering{\mbox{\psfig{figure=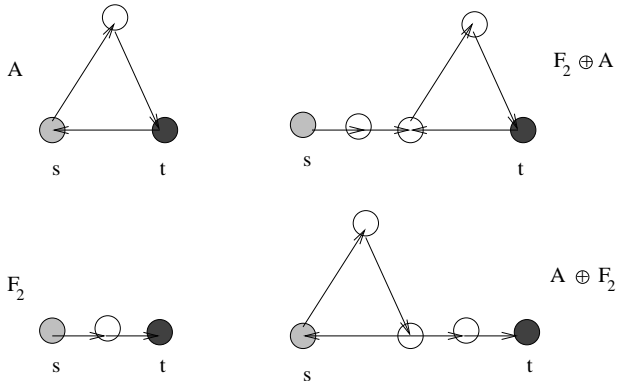}}}
\caption{Example showing the non-commutativity of addition in ${\mathcal F}$.}
\label{non-comm-plus-figure}
\end{figure}
\end{example}

The previous example proves the next lemma.
\begin{lemma}
\index{$\oplus$!non-commutativity}
The operation $\oplus$ is not commutative.
\end{lemma}

% MOD2
\begin{definition}[$\oplus$-Irreducible]
\label{oplus irreducibile} 
\index{$\oplus$!irreducibile}
\index{irreducible!$\oplus$} 
A flow graph $A$ is called {\em $\oplus$-reducible} if there is
decomposition of $A$ as sum of non-trivial flow graphs $B$ and $C$.
Otherwise $A$ is called $\oplus$-irreducible.
\end{definition}

% MOD3
\begin{definition}[$s$-standard and $t$-standard]
\label{standardness}
\index{$s$ standard} 
\index{$t$ standard} 
\index{standard!$s,t$} 
Let $A$ be a flow graph.  $A$ is called {\em $s$-standard} if it can
not be decomposed as $A = B \oplus C$, where $B$ is infinitesimal, and
$C$ is an arbitrary (possibly trivial) flow graph.  $A$ is called 
% REV IM
{\em $t$-standard} if it cannot be decomposed as $A = B \oplus C$,
where $B$ is an arbitrary (possibly trivial) flow graph, and and $C$
is an infinitesimal flow graph.
\end{definition}

% MOD3
Note that $\oplus$-irreducible non-infinitesimal flow graphs are both
$s$-standard and $t$-standard.  On the other hand,
$\oplus$-irreducible infinitesimal flow graphs are neither
$s$-standard nor $t$-standard.  There are no restrictions on the
$s$-standardness and $t$-standardness properties of general
$\oplus$-reducible flow graphs.

% MOD2
We would like to devise a graph-theoretic characterization of
$\oplus$-irreducibility.  Towards this, the next definition is the
flow graph analogue of a cut vertex in standard graphs.

\begin{definition}[Splitting vertex for a flow graph]
% MOD3
\label{splitting vertex for a flow graph}
\index{splitting!vertex for a flow graph} 
\index{vertex!splitting---for a flow graph} 
\index{flow graph!vertex splitting for a ---} 
We say that $w$ is a {\em splitting vertex} for flow graph $A=(G_A,
s_A, t_A)$ if $w \neq s_A,t_A$ and the deletion of $w$ from $G_A$
produces {\em at least two} non-trivial components, with $s_A$ and $t_A$ lying in
distinct components.  We denote the component containing $s_A$ as
$G_A^s(w)$, the one containing $t_A$ as $G_A^t(w)$, and the remaining
components as $G_A^{\epsilon}(w)$.  Note that $G_A^{\epsilon}(w)$ may
be the union of several disjoint components, and hence is not
necessarily connected.  Let $i_s:G_A^s(w) \hookrightarrow G_A$,
$i_t:G_A^t(w) \hookrightarrow G_A$, and $i_\epsilon:G_A^\epsilon(w)
\hookrightarrow G_A$ denote the natural subgraph injections.  Since
$w$ is a splitting vertex for $A$, $s_A \in$ Im$(i_s)$, $t_A \in$
Im$(i_t)$.  The pair of graphs $(G_A^s(w),G_A^t(w))$ is called the
$(s_A,t_A)$ splitting of $G_A$ {\em induced} by $w$.
% MODIFICATIONS NECESSARY: SEARCH FOR 'splitting' TO FIND SIDE EFFECTS
\end{definition}
% MOD3
\begin{remark}
\label{split-cut-different}
If $w$ is a splitting vertex in flow graph $A=(G_A, s_A, t_A)$,
then by Definition \ref{splitting vertex for a flow graph}, $w$ is a
cut vertex in $G_A$.  The converse is false, however, since not every
cut vertex in $G_A$ is a splitting vertex in $A$.  

For a concrete example, the reader may wish to consider the flow graph
$A$ in Figure~\ref{plus-general-picture} on
page~\pageref{plus-general-picture}, where the radius $1$ sphere of
$t_A$ contains two cut vertices for $G_A$, only one of which is a
splitting vertex for $A$.

As a more general example, consider an infinitesimal flow graph
$A=(G_A,s_A,t_A)$, for which $G_A$ is 1-connected (as a graph).  Since
$A$ is infinitesimal, $s_A=t_A$, so $A$ does not possess a splitting
vertex.  By 1-connectedness, however, $G_A$ contains a cut vertex.
\end{remark}

\begin{definition}[Flow graph splitting]
\label{flow graph splitting}
\index{flow graph!splitting}
\index{splitting!flow graph}
% MOD3
Suppose vertex $w$ is a splitting vertex for flow graph $A=(G_A, s_A,
t_A)$.  Take $t_{A_s^w}$ to be a new vertex (not present in
$V[G_A^s(w)] \cup V[G_A^\epsilon(w)]$), and define flow graph $A_s^w =
(G_{A_s^w}, s_{A_s^w}, t_{A_s^w})$ as follows:
\begin{eqnarray*}
V[G_{A_s^w}] & = & V[G_A^s(w)] \cup V[G_A^\epsilon(w)] \cup \{t_{A_s^w}\},\\
E[G_{A_s^w}] & = & E[G_A^s(w)] \cup E[G_A^\epsilon(w)] \\
             &   & \cup \; \{(t_{A_s^w},u) \;|\; (w,u) \in E[G_A], u \in
(V[G_A^s(w)] \cup V[G_A^\epsilon(w)]) \} \\
             &   & \cup \; \{(u,t_{A_s^w}) \;|\; (u,w) \in E[G_A], u \in
(V[G_A^s(w)] \cup V[G_A^\epsilon(w)]) \},\\
s_{A_s^w}    & = & i_s^{-1}(s_A).
\end{eqnarray*}
Analogously, let $s_{A_t^w}$ be a vertex not present in
$V[G_A^t(w)]$.  Define $A_t^w$ to be the flow graph
$(G_{A_t^w}, s_{A_t^w}, t_{A_t^w})$ as follows.
\begin{eqnarray*}
V[G_{A_t^w}] & = & V[G_A^t(w)] \cup \{s_{A_t^w}\},\\
E[G_{A_t^w}] & = & E[G_A^t(w)] \\
             &   & \cup \; \{(s_{A_t^w},u) \;|\; (w,u) \in E[G_A], u \in
V[G_A^t(w)]\} \\
             &   & \cup \; \{(u,s_{A_t^w}) \;|\; (u,w) \in E[G_A], u \in
V[G_A^t(w)]\},\\
t_{A_t^w}    & = & i_t^{-1}(t_A).
\end{eqnarray*}
The pair of flow graphs $(A_s^w,A_t^w)$ is referred to as the {\em
splitting} of $A$ induced by $w$.
\end{definition}

\begin{lemma}[$\oplus$-Irreducibility Lemma]
\label{oplus Irreducibility Lemma} 
\index{$\oplus$!Irreducibility Lemma} 
\index{Irreducibility Lemma!$\oplus$} 
Flow graph $A=(G_A, s_A, t_A)$ is $\oplus$-reducible if and only if
$V[G_A]$ contains a splitting vertex for $A$.
\end{lemma}
\begin{proof}
If $A = B \oplus C$, then
\begin{eqnarray*}
w = {\sigma^{\NStinyoplus-1}_{t_B \approx s_C}} (t_B) =
{\tau^{\NStinyoplus-1}_{t_B \approx s_C}} (s_C)
\end{eqnarray*}
is a splitting vertex for $A$ (see expression~(\ref{plus-injections})
on page \pageref{plus-injections} for definitions of the $\sigma$ and
$\tau$ injections).  Conversely, if $w$ is a splitting vertex for $A$,
then $A = A_s^w \oplus A_t^w$.
\end{proof}

% MOD2
Suppose $A = A_0 \oplus A_1$ and $B = B_0 \oplus B_1$ are two flow
graphs. If $A_i=B_i$ ($i=0,1$) then $A=B$.  We would like to
investigate the extent to which the converse is true.  Towards this,
we introduce the following property.

% MOD2
\begin{definition}(Undirected Form)
\label{undirected form}
Let $G=(V,E)$ be a directed multigraph.  Then the {\em undirected
form} $U(G)=(V, \bar{E})$ is defined to be the undirected multigraph
on $V$, in which there is an undirected edge $(u,v)$ in $\bar{E}$ for
each directed edge from $u$ to $v$ or from $v$ to $u$ in $E$.  
% MOD3
For any directed edge $e \in E$, its corresponding {\em undirected}
representative in $\bar{E}$ is denoted as $U(e)$. In short, $U$ is an
operation on directed edges which ``forgets'' their orientations. We
extend $U$ to act on sequences of edges from $E$ in the obvious
manner.
% MOD3
For any undirected edge $e' \in \bar{E}$, its corresponding {\em
directed} representative in $E$ is denoted as $D_G(e')$.  In short,
$D_G$ is an operation on undirected edges which assigns orientations
according to the orientations of edges in $G$. We extend $D_G$ to act
on sequences of edges from $\bar{E}$ in the obvious manner.
\end{definition}

% MOD2
\begin{definition}
\label{st-property}
\index{{\sc st} Property}
\index{{\sc st}-flow graph}
\index{$F^{ST}$}
\index{flow graph!{\sc st} Property}
Let $A=(G_A, s_A, t_A)$ be an arbitrary flow graph.  An edge $e$ in
$E[G_A]$ is said to have the {\sc st} Property if in the undirected graph
$U(G_A)$ there is a (non self-intersecting) {\em path} from $s_A$ to
$t_A$ which traverses $U(e)$. We say $A$ is an {\sc st}-flow graph if
every edge $e$ in $E[G_A]$ has the {\sc st} Property.  The set of all
{\sc st}-flow graphs will be denoted $F^{ST} \subset F$.
\end{definition}

% MOD3
Clearly, no {\sc st}-flow graph is infinitesimal.  Hence if the source
and target vertices of an {\sc st}-flow graph coincide, then it must
be trivial.  {\sc st}-flow graphs have several nice features which
will be useful.

% MOD2
\begin{lemma}
\label{st-closure}
$F^{ST}$ is closed under $\oplus$.
\end{lemma}
\begin{proof}
If $A$ and $B$ are {\sc st}-flow graphs, then any edge $e$ in $G_{A
\tinyoplus B}$ is either in $Im({\sigma^{\NStinyoplus}_{t_A \approx
s_B}})$ or in $Im({\tau^{\NStinyoplus}_{t_A \approx s_B}})$.  Suppose
$e$ is in $Im({\sigma^{\NStinyoplus}_{t_A \approx s_B}})$.  Since $A$
is an {\sc st}-flow graph, there is a path $p$ from $s_A$ to $t_A$ in
$U(G_A)$ which contains $U({\sigma^{\NStinyoplus-1}_{t_A \approx
s_B}}(e))$; and since $B$ is an {\sc st}-flow graph, there is a path
$q$ in $U(G_B)$ from $s_B$ to $t_B$.  Then $p$ concatenated with $q$
is a path in $U(G_{A \tinyoplus B})$ which contains $U(e)$ and
connects $s_{A \tinyoplus B}$ to $t_{A \tinyoplus B}$.  The case when
$e$ is in $Im({\tau^{\NStinyoplus}_{t_A \approx s_B}})$ is
analogous.
\end{proof}

The next lemma considers the reverse implication.

\begin{lemma}
\label{st-closure-reverse}
Let $A=(G_A,s_A,t_A)$ and $B=(G_A,s_A,t_A)$ be flow graphs.  If $A
\oplus B$ is an {\sc st}-flow graph, then $A$ and $B$ are {\sc
st}-flow graphs.
\end{lemma}
\begin{proof} 
Suppose, towards contradiction, that $A$ is not an {\sc st}-flow
graph.  Then there is an edge $e$ that does not have the {\sc st}
Property in $A$.  So $e$ is not on any path from $s_A$ to $t_A$ in
$G_A$.

Suppose, towards contradiction, that the edge
$\sigma^{\NStinyoplus}_{t_A \approx s_B}(e)$ is on a path $p$ from
$s_{A\tinyoplus B}$ to $t_{A\tinyoplus B}$ in $G_{A\tinyoplus B}$.
Then, the restriction of this $p$ to $Im(\sigma^{\NStinyoplus})$
contradicts the fact that $e$ does not have the {\sc st} Property.  So
there is no path from $s_{A\tinyoplus B}$ to $t_{A\tinyoplus B}$ in
$G_{A\tinyoplus B}$ containing $\sigma^{\NStinyoplus}_{t_A
\approx s_B}(e)$.  

Thus $\sigma^{\NStinyoplus}_{t_A \approx s_B}(e)$ does not have the
{\sc st} Property in $A \oplus B$.  This is a contradiction, hence $A$
must be an {\sc st}-flow graph.  The case where we assume $B$ is not
an {\sc st}-flow graph leads to a similar contradiction.
\end{proof}

% REV IM
Let Let $A=(G_A,s_A,t_A)$ be an arbitrary flow graph.  Consider two
edges $e$ and $e'$ in $G_A$.  Suppose that $e$ does not have the {\sc
st} Property.  Then $e$ cannot lie inside any path from $s$ to $t$ in
$U(G)$.
\begin{itemize}
\item Suppose $e'$ has the {\sc st} Property in $G$.  Then there is path $p$
from $s$ to $t$ in $U(G)$ containing $e'$.  Certainly $p$ cannot
contain $e$, since $e$ does not have the {\sc st} Property.  So $p$
survives in $G\backslash e$, and hence $e'$ has the {\sc st} Property
in $G\backslash e$.
\item Suppose $e'$ does not have the {\sc st} Property in $G$.  Then there is
no path $p$ from $s$ to $t$ in $U(G)$ containing $e'$.  But then the
deletion of $e$ from $G$ does not change this fact, since deletion of
edges cannot create new paths.  It follows that in $e'$ does not have
the {\sc st} Property in $G\backslash e$.
\end{itemize}
We have shown that if an edge $e$ does not have the {\sc st} Property,
then its deletion does not affect other edges $e'$, for whom the {\sc
st} Property stays unchanged from $G$ to $G\backslash e$.  It follows
that sequentially deleting all edges which violate the {\sc st}
Property of $G$ yields a unique directed connected graph whose vertex
set contains the source and target vertices of the original flow
graph.  This leads to the next definition, by which we can approximate
general flow graphs using suitable large {\sc st}-flow subgraphs.

% MOD2
\begin{definition}
\label{st-core}
\index{{\sc st}-Core}
Let $A=(G_A, s_A, t_A)$ be a general flow graph.  We define the {\em
{\sc st}-core} of $A$ to be the largest {\sc st}-flow graph
$\hat{A}=(\hat{G}_A, s_A, t_A)$ contained in $A$, where $\hat{G}_A =
(\hat{V}[G_A], \hat{E}[G_A])$ is obtained by deleting all edges from
$G_A$ which do not have the {\sc st} Property and considering the
connected subgraph induced by the remaining edges.  Note that
$\hat{A}$ is an {\sc st}-flow graph on a vertex set $\hat{V}[G_A]
\supseteq \{s_A, t_A\}$.
\end{definition}

% MOD2
\begin{lemma}
\label{core-distributes}
The {\sc st}-core operation distributes over $\oplus$.  
\end{lemma}
\begin{proof}
We would like to show that $\hat{A} \oplus \hat{B}$ is isomorphic to
the {\sc st}-core($A\oplus B$).

Let $e$ be an arbitrary edge in the {\sc st}-core($A\oplus B$).  Then
there is a path $p$ from $s_{A\tinyoplus B}$ to $t_{A\tinyoplus B}$ in
$U(G_{A\tinyoplus B})$, such that $p$ contains $U(e)$.  But $p$ must
pass through the cut vertex $w=\sigma^{\NStinyoplus}_{t_A \approx
s_B}(t_A)$ in $G_{A\tinyoplus B}$.
\begin{itemize}
\item[(i)] If $e$ is in $Im(\sigma^{\NStinyoplus}_{t_A \approx s_B})$,
let $p'$ be the initial segment of $p$ which connects $s_{A\tinyoplus
B}$ to $w$ in $U(G_{A\tinyoplus B})$; then
$U(\sigma^{\NStinyoplus-1}_{t_A \approx s_B}(D_{G_{A\tinyoplus B}}(p')))$ is a path
containing $U(e)$ connecting $s_A$ to $t_A$ in $U(G_A)$, so
$e'=\sigma^{\NStinyoplus-1}_{t_A \approx s_B}(e)$ is in $\hat{A}$.

\item[(ii)] If $e$ is in $Im(\tau^{\NStinyoplus}_{t_A \approx s_B})$,
let $p'$ be the final segment of $p$ which connects $w$ to
$t_{A\tinyoplus B}$ in $U(G_{A\tinyoplus B})$; then $U(\tau^{\NStinyoplus-1}_{t_A \approx
s_B}(D_{G_{A\tinyoplus B}}(p')))$ is a path containing $U(e)$ connecting $s_B$ to $t_B$ in $U(G_B)$, so
$e'=\tau^{\NStinyoplus-1}_{t_A \approx s_B}(e)$ is in $\hat{B}$.
\end{itemize}
This shows that $\hat{A} \oplus \hat{B}$ contains {\sc
st}-core($A\oplus B$).

Conversely, let $e'$ be an edge in $\hat{A}$.  Then there is a path
$p'$ containing $e'$ which connects $s_A$ to $t_A$ in $U(G_A)$.  It
follows that $U(\sigma^{\NStinyoplus}_{t_A \approx s_B}(D_{G_A}(p')))$
connects $s_{A\tinyoplus B}$ to the cut vertex
$w=\sigma^{\NStinyoplus}_{t_A \approx s_B}(t_A)$ in $U(G_{A\tinyoplus
B})$.  Let $p''$ be a path from $s_B$ to $t_B$ in $U(G_B)$.  Then by
concatenating $U(\sigma^{\NStinyoplus}_{t_A \approx s_B}(D_{G_A}(p')))$ with
$U(\tau^{\NStinyoplus}_{t_A \approx s_B}(D_{G_B}(p'')))$ we obtain a path
$p$ connecting $s_{A\tinyoplus B}$ to $t_{A\tinyoplus B}$ through
$\sigma^{\NStinyoplus}_{t_A \approx s_B}(e)$ in $U(G_{A\tinyoplus
B})$.  A similar argument can be carried out when $e'$ is an edge in
$\hat{B}$.  This shows that the {\sc st}-core($A\oplus B$) contains
$\hat{A} \oplus \hat{B}$.
\end{proof}

The next lemma stands in contrast to Remark \ref{split-cut-different}.

% MOD2
\begin{lemma}
\label{spanning-cut-same}
Let $A=(G_A, s_A, t_A)$ be an {\sc st}-flow graph, and let $w \neq
s_A, t_A$ be a cut vertex of $G_A$.  Then $w$ is a splitting vertex of
$A$.
\end{lemma}
\begin{proof}
Let $w$ be a cut vertex of $G_A$.  Delete $w$ and denote the union of
those components that contain neither $s_A$ nor $t_A$ as
$G_A^{\epsilon}(w)$ (see Definition \ref{splitting vertex for a flow
graph}).  Towards contradiction, suppose $V[G_A^{\epsilon}(w)]\neq
\emptyset$; let $v$ be any vertex therein.  Then any walk from $s_A$
to $t_A$ in $U(G_A)$ through $v$ must visit $w$ twice and hence
self-intersect.  It follows that edges incident to $w$ violate the
assumption that $A$ has Property {\sc st} (see Definition
\ref{st-property}).  Hence it must be that $V[G_A^{\epsilon}(w)] =
\emptyset$.

Since $w$ is a cut vertex of $G_A$ it's deletion generates at least
two components.  Towards contradiction, let $s_A$ and $t_A$ lie in the
same component.  Then there is at least one component which contains
neither $s_A$ nor $t_A$, and so $V[G_A^{\epsilon}(w)] \neq \emptyset$.
We have shown that the deletion of $w$ generates precisely two
components, and that $s_A$ and $t_A$ lie in distinct components.
Hence $w$ is a splitting vertex of $A$.
\end{proof}
% MOD2
\begin{lemma}[Crossing Summands Lemma]
\label{spanning-lemma}
Let $A=(G_A, s_A, t_A)$ and $B=(G_B, s_B, t_B)$ be flow graphs.  Let
$C=(G_C, s_C, t_C)$ be an {\sc st}-flow graph with
\begin{eqnarray*}
   \alpha:G_C \hookrightarrow G_{A \tinyoplus B}
\end{eqnarray*}
be a graph embedding satisfying $\alpha(s_C) = s_{A\tinyoplus B}$ and $\alpha(G_C)
\not\subset Im({\sigma^{\NStinyoplus}_{t_A \approx s_B}})$.  Then
\begin{eqnarray*}
\alpha^{-1}\sigma^{\NStinyoplus}_{t_A \approx s_B}(t_A)=\alpha^{-1}\tau^{\NStinyoplus}_{t_A \approx s_B}(s_B)
\end{eqnarray*}
is a splitting vertex for $C$.
\end{lemma}
\begin{proof}
Let $w = \sigma^{\NStinyoplus}_{t_A \approx
s_B}(t_A)=\tau^{\NStinyoplus}_{t_A \approx s_B}(s_B)$.  Clearly, $w$
is a splitting vertex in $A \oplus B$ and hence a cut vertex in $G_{A
\tinyoplus B}$.  Let $v$ be a vertex in $\alpha(G_C)$ that is not in
$Im({\sigma^{\NStinyoplus}_{t_A \approx s_B}})$.  Since $w$ is a cut
vertex in $G_{A \tinyoplus B}$, every path from $s_{A \tinyoplus B}$
to $v$ in $U(G_{A \tinyoplus B})$ must pass through $w$.  Given such a
path
\begin{eqnarray*}
p=(p_0=s_{A \tinyoplus B}, p_1, \ldots p_{l-1} , p_l = w, p_{l+1}, \ldots, p_{|p|-1}, p_{|p|} = v),
\end{eqnarray*}
which lies entirely in $U(Im(\alpha))$, we define the pullback
$\alpha^{-1}(p)$ to be a path in $U(G_C)$ in the obvious manner.  By
considering pullbacks of paths from $s_{A \tinyoplus B}$ to $v$ in
$U(G_{A \tinyoplus B})$, we see that every path from $s_{C}$ to
$\alpha^{-1}(v)$ in $U(G_C)$ must pass through $\alpha^{-1}(w)$.
Hence the deletion of $\alpha^{-1}(w)$ must generate at least two
components, one of which contains $s_{C}$ while another contains
$\alpha^{-1}(v)$.  This shows that $\alpha^{-1}(w)$ is a cut vertex in
$U(G_C)$.  Since $C$ is an {\sc st}-flow graph, Lemma
\ref{spanning-cut-same} holds and $\alpha^{-1}(w)$ is a splitting
vertex for $C$.
\end{proof}

\begin{figure}[bth]
\centering{\mbox{\psfig{figure=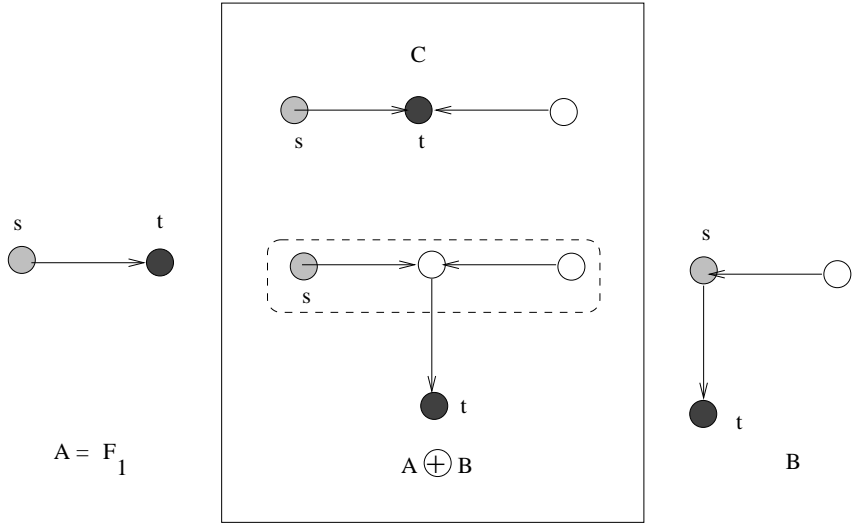}}}
\caption{An example where the Crossing Summands Lemma fails for general flow graphs.}
\label{ssl-fails}
\end{figure}
% REV IM
Consider the flow graphs $A$, $B$, $C$ shown in Figure
\ref{ssl-fails}.  Note that $C$ is not an {\sc st}-flow graph, and
$\alpha:G_C \hookrightarrow G_{A \tinyoplus B}$ embeds as shown.
Since $t_C$ is not a splitting vertex for $C$, the Crossing Summands
Lemma is seen to fail when $C$ is a general flow graph.

\begin{proposition}[Component-wise decomposition of isomorphisms under
$\oplus$]
\label{component-wise-decomposition-of-isomorphisms-under-plus}
\index{component-wise decomposition!of isomorphisms}
\index{isomorphisms!component-wise}
Suppose $A$ and $B$ flow graphs, expressed as sums of
$\oplus$-irreducible {\sc st}-flow graphs as follows:
\begin{eqnarray*}
A & = & A_0 \oplus A_1 \oplus \cdots \oplus A_{m-1},\\
B & = & B_0 \oplus B_1 \oplus \cdots \oplus B_{n-1}.
\end{eqnarray*}
Then
\begin{eqnarray*}
A = B & \Leftrightarrow & m=n \text{ and } A_i = B_i \text{ for all $i
= 0, \ldots , m-1$.}
\end{eqnarray*}
\end{proposition}
\begin{proof}
[$\Leftarrow$] Let $\phi_i:A_i \rightarrow B_i$ be given
component-wise isomorphisms, for $i=0, \ldots , m-1$.  Define $\phi: A
\rightarrow B$ by defining $\phi|_{A_i} = \phi_i$.  Since $t_{B_i} =
\phi_i(t_{A_i}) = \phi_{i+1}(s_{A_{i+1}}) = s_{B_{i+1}}$ for $i=0,
\ldots, m-2$, this provides a well-defined isomorphism between $A$ and
$B$.  Indeed, this implication holds for arbitrary flow graphs.

[$\Rightarrow$] We prove the statement by induction on $\max(m,n)$.
In the case when $m=n=1$, the claim is trivial.  For the inductive
step, let $\phi: A \rightarrow B$ be an isomorphism.  Consider
$\phi(A_0)$, and take $k$ to be the smallest integer in $\{0, \ldots,
n-1\}$ for which $\phi(A_0)$ is a subgraph of $B_0 \oplus B_1 \oplus
\cdots \oplus B_k$.  Since $A_0$ is $\oplus$-irreducible it contains
no splitting vertices.  
% MODIFICATIONS NECESSARY
Since $A_0$ is an {\sc st}-flow graph, it must be that $k=0$;
otherwise by Lemma \ref{spanning-lemma} $A_0$ would contain a
splitting vertex, and so (by Lemma \ref{oplus Irreducibility Lemma})
be $\oplus$-reducible.  Since $k=0$, we have shown that $\phi(A_0)$ is
a subgraph of $B_0$.  Now, repeating the argument for $B_0$ using
$\phi^{-1}$, we see that $\phi^{-1}(B_0)$ is a subgraph of $A_0$.  It
follows that $A_0$ is isomorphic to $B_0$ under a suitable restriction
of $\phi$.  Now, since $\phi(A) = B$ and $\phi(A_0) = B_0$ it follows
that $\phi(A \backslash A_0) = B \backslash B_0$, or more specifically
$\phi(A_1 \oplus \cdots \oplus A_{m-1}) = B_1 \oplus \cdots \oplus
B_{n-1}$.  By inductive hypothesis, this implies that $m=n$ and $A_i =
B_i$ for all $i=0, \ldots, m-1$.
\end{proof}

% REV IM
Because the Crossing Summands Lemma fails (and lies at the heart of
the proof of Proposition
\ref{component-wise-decomposition-of-isomorphisms-under-plus}) the
latter Proposition also fails for general flow graphs.  The reader can
verify for example, from Figure \ref{ssl-fails}, that $A \oplus B = C
\oplus A$ but $A\neq C$ and $B \neq A$.

\begin{definition}[Splitting vertex ranking]
\label{splitting vertex ranking}
\index{splitting! vertex ranking}
Given flow graph $A=(G_A, s_A, t_A)$, let $\chi(A) \subset V[G_A]$ be
the set of all splitting vertices for $A$.  We define the {\em
$s$-ranking} and {\em $t$-ranking} functions $r_s^A, r_t^A: \chi(A)
\rightarrow {\mathbb N}$ as follows:
\begin{eqnarray*}
r_s^A(w) & = & | V[G_A^s(w)] \cap \chi(A) |,\\
r_t^A(w) & = & | V[G_A^t(w)] \cap \chi(A) |.
\end{eqnarray*}
When it is clear from the context, we denote $r_s(w) =
r_s^A(w)$ and $r_t(w) = r_t^A(w)$.
\end{definition}

\begin{lemma}
\label{relative-ranks-of-split-vertices}
Let $w \in \chi(A)$  be a splitting vertex for flow graph $A=(G_A, s_A,
t_A)$.
Then for all $u \in V[G_A^s(w)] \cap \chi(A)$:
\begin{eqnarray*}
r_s(u) & < & r_s(w),\\
r_t(u) & > & r_t(w);
\end{eqnarray*}
and for all $u \in V[G_A^t(w)] \cap \chi(A)$:
\begin{eqnarray*}
r_s(u) & > & r_s(w),\\
r_t(u) & < & r_t(w).
\end{eqnarray*}
\end{lemma}
\begin{proof}
First, note that for any $u$ in $(V[G_A^s(w)] \cup V[G_A^t(w)]) \cap
\chi(A)$
\begin{eqnarray}
\label{total rank}
r_s(u) + r_t(u) + 1 & = & |\chi(A)|
\end{eqnarray}
Now if $u \in V[G_A^s(w)] \cap \chi(A)$, then since $w \in (V[G_A^t(u)]
\backslash
V[G_A^t(w)]) \cap \chi(A)$, so it follows that $V[G_A^t(w)]
\subsetneq V[G_A^t(u)]$.  But then $r_t(w) < r_t(u)$.  By
expression~\ref{total rank} above, it follows that $r_s(w) > r_s(u)$.
The proof for the case when $u \in V[G_A^t(w)] \cap \chi(A)$ is
analogous.
\end{proof}

\begin{lemma}
\label{ordering-scheme-for-splitting-vertices}
Given a flow graph $A=(G_A, s_A, t_A)$, for each $i = 0, 1, \ldots,
|\chi(A)|-1$ there is a unique vertex $v_i$ in $\chi(A)$ with the
property that $r_s(v_i) = i$.
\end{lemma}
\begin{proof}
First we note that one cannot have two distinct vertices $v$, $v'$
having $r_s(v) = r_s(v')$, since either $v \in V[G_A^{s}(v')]$ or $v'
\in V[G_A^{s}(v)]$, and so by
Lemma~\ref{relative-ranks-of-split-vertices} it follows that $r_s(v)
\neq
r_s(v')$.
Base case: $i=0$.  Let $w_0$ be any vertex in $\chi(A)$. If $r_s(w) >
0$, then $V[G_A^s(w)] \cap \chi(A)$ is not empty.  So let $w_1$ be any
vertex in $V[G_A^s(w)] \cap \chi(A)$.  By
Lemma~\ref{relative-ranks-of-split-vertices}, $r_s(w_1) < r_s(w_0)$.
Repeating in this fashion, after finitely many steps $w_0
\rightsquigarrow w_1 \rightsquigarrow \ldots $ we find some vertex
$v_0$ for which $r_s(v_0) = 0$.
Inductive step $i+1$: Let $v_i$ be the unique vertex in $\chi(A)$
having $r_s(v_i) = i$.  Define $v_{i+1}$ to be the vertex in
$V[G_A^{t}(v_i)] \cap \chi(A)$ for whose $s$-rank is minimal.  Since
\begin{eqnarray*}
V[G_A^{s}(v_{i+1})] \cap \chi(A) & = & [V[G_A^{s}(v_{i})] \cap \chi(A)]
\cup \{ v_i \},
\end{eqnarray*}
it follows that $r_s(v_{i+1}) = r_s(v_i)+1 = i+1$, hence the result.
\end{proof}

\begin{definition}[Canonical $\oplus$-decomposition]
\label{canonical-oplus-decomposition}
\index{canonical $\oplus$-decomposition}
\index{$\oplus$-decomposition}
\index{decomposition!$\oplus$}
Let $A=(G_A, s_A, t_A)$ be a flow graph.  Take $\chi(A) = \{v_0, v_1,
\ldots, v_{|\chi(A)|-1}\}$ to be the set of splitting vertices for
$A$, ordered according to the indexing scheme postulated in
Lemma~\ref{ordering-scheme-for-splitting-vertices}.  Define $A^{(0)} =
A_s^{v_0}$, $\bar{A}^{(0)} = A_t^{v_0}$, and then for each $i=1, 2,
\ldots, |\chi(A)|-1$, put
\begin{eqnarray*}
      A^{(i)} & = & {(\bar{A}^{(i-1)})}_s^{v_{i}},\\
\bar{A}^{(i)} & = & {(\bar{A}^{(i-1)})}_t^{v_{i}}.
\end{eqnarray*}
We shall denote $\bar{A}^{(|\chi(A)| - 1)}$ as $A^{(|\chi(A)|)}$.  The
{\em canonical $\oplus$-decomposition of $A$} is defined to be the
sequence
\begin{eqnarray*}
   \langle A \rangle \definemathit (A^{(0)}, A^{(1)}, \ldots,
A^{(|\chi(A)|-1)}, A^{(|\chi(A)|)}).
\end{eqnarray*}
Note that the effectiveness of this definition guarantees uniqueness
of the decomposition.
\end{definition}

% MODIFICATIONS NECESSARY
\begin{lemma}[$\oplus$-decompositions for sums]
\index{canonical $\oplus$-decomposition!for sums}
\label{canonical-plus-decomp-for-sum}
Given $A=(G_A, s_A, t_A)$ and $B=(G_B, s_B, t_B)$, two flow graphs
with their respective canonical $\oplus$-decompositions $\langle A
\rangle$ and $\langle B \rangle$. If $B$ is $s$-standard, then the
canonical $\oplus$-decomposition of $A \oplus B$ is $\langle A
\rangle\langle B \rangle$, the concatenation of $\langle A \rangle$
with $\langle B \rangle$.
\end{lemma}
\begin{proof}
First, note that $|\chi(A \oplus B)| = |\chi(A)| + |\chi(B)| + 1$.  More
specifically, if
\begin{eqnarray*}
\chi(A) & = & \{u_0, u_1, \ldots, u_{|\chi(A)|-1}\}, \text{ and}\\
\chi(B) & = & \{v_0, v_1, \ldots, v_{|\chi(B)|-1}\}
\end{eqnarray*}
are the sets of splitting vertices for $A$ and $B$ respectively,
ordered by ascending $s$-rank, according to the indexing scheme
postulated in Lemma~\ref{ordering-scheme-for-splitting-vertices}, then
$A \oplus B$ has splitting vertices:
\begin{eqnarray*}
\chi(A \oplus B) & = & \left\{ {\sigma^{\NStinyoplus}_{t_A \approx s_B}}
(u_0), {\sigma^{\NStinyoplus}_{t_A \approx s_B}} (u_1), \ldots,
{\sigma^{\NStinyoplus}_{t_A \approx s_B}} (u_{|\chi(A)|-1}),\right. \\
& & \phantom{\{} \; {\sigma^{\NStinyoplus}_{t_A \approx s_B}} (t_A) =
{\tau^{\NStinyoplus}_{t_A \approx s_B}} (s_B),\\
& & \phantom{\{} \left. {\tau^{\NStinyoplus}_{t_A \approx s_B}} (v_0),
{\tau^{\NStinyoplus}_{t_A \approx s_B}} (v_1), \ldots, {\tau^{\NStinyoplus}_{t_A
\approx s_B}} (v_{|\chi(B)|-1}) \right\}.\\
\end{eqnarray*}
But since ${\sigma^{\NStinyoplus}_{t_A \approx s_B}}$ and
${\tau^{\NStinyoplus}_{t_A \approx s_B}}$ are injections,
\begin{eqnarray*}
(A\oplus B)^{(i)} =
\left\{
\begin{array}{ll}
A^{(i)}                   & \text{for } 0 \leqslant i < |\chi(A)|,\\
\bar{A}^{(|\chi(A)|-1)}   & \text{for } i = |\chi(A)|,\\  
B^{(i-|\chi(A)|-1)}       & \text{for } |\chi(A)| < i \leqslant
|\chi(A)|+|\chi(B)|,\\
\bar{B}^{(|\chi(B)|-1)}   & \text{for } i = |\chi(A)|+|\chi(B)|+1.
\end{array}
\right.
\end{eqnarray*}
% MODIFICATIONS NECESSARY
It follows that $\langle A \oplus B \rangle = \langle A
\rangle\langle B \rangle$.
\end{proof}

% REV IM
We note that if $B$ is not $s$-standard then $(A\oplus
B)^{(|\chi(A)|)}$ is not equal to $\bar{A}^{(|\chi(A)|-1)}$ and the
Lemma \ref{canonical-plus-decomp-for-sum} fails to hold.  The reader
can verify this by considering the flow graphs $A$ and $B$ shown in
Figure \ref{ssl-fails}.

\begin{proposition}[Correctness of the $\oplus$-decomposition]
\label{decomposition-is-valid}
Consider the canonical $\oplus$-decomposition of
$A$ as given in Definition~\ref{canonical-oplus-decomposition}:
\begin{eqnarray*}
\langle A^{(0)}, A^{(1)}, \ldots, A^{(|\chi(A)|)} \rangle.
\end{eqnarray*}
Then $A = A^{(0)} \oplus A^{(1)} \oplus A^{(2)} \cdots  \oplus A^{(|\chi(A)|-1)} \oplus A^{(|\chi(A)|)}$, 
and every summand is $\oplus$-irreducible.
\end{proposition}
\begin{proof}
Since $r_s(v_i) = i$, it follows that $\chi(A^{(i)}) = \emptyset$
for all $i = 0, 1, \ldots, |\chi(A)|$.  Since each summand has no
splitting vertices, by Lemma~\ref{oplus Irreducibility Lemma},
each is $\oplus$-irreducible.
We prove the Proposition by induction on $|\chi(A)|$.
The base case when $|\chi(A)|=1$ is straightforward, since
Definition~\ref{canonical-oplus-decomposition} specified $A^{(0)} =
A_s^{v_0}$
and $\bar{A}^{(0)} = A_t^{v_0}$.  Then, since $A_s^{v_0} \oplus
A_t^{v_0} = A$ for any splitting vertex $v_0$, the result follows.
Suppose the Proposition has been proved for all flow graphs $B$ which
enjoy $|\chi(B)| \leqslant k$.  Let $A$ be a flow graph with
$|\chi(A)| = k+1$.  Unravelling
Definition~\ref{canonical-oplus-decomposition} yields $A^{(k+1)} =
\bar{A}^{(k)}
= (\overline{A}^{(k-1)})_t^{v_{k}} = A_t^{v_{k}}$.  Then since $A =
A_s^{v_{k}} \oplus A_t^{v_{k}}$ and $|\chi(A_s^{v_{k}})| = k$, by
inductive hypothesis and Lemma~\ref{canonical-plus-decomp-for-sum},
\begin{eqnarray*}
A & = & A_s^{v_{k}} \oplus A_t^{v_{k}}\\
  & = & \langle A_s^{v_{k}} \rangle \oplus A_t^{v_{k}} \\
  & = & (A_s^{v_{k}})^{(0)} \oplus \ldots (A_s^{v_{k}})^{(k)} \oplus
A_t^{v_{k}}\\
  & = & A^{(0)} \oplus A^{(1)} \oplus A^{(2)} \cdots  \oplus A^{(k)}
\oplus A^{(k+1)}.
\end{eqnarray*}
The result follows.
\end{proof}

% REV IM
% MODIFICATIONS NECESSARY
% \begin{lemma}({\sc st}-core Extention)
% \label{extension}
% Let $\psi$ be flow graph isomorphism of $A=(G_A, s_A, t_A)$ and
% $\hat{\phi}: \hat{A} \rightarrow \hat{A}$ be a flow graph isomorphism
% of it's {\sc st}-core.  Then one can extend $\hat{\phi}$ to an
% isomorphism $\phi : A \rightarrow A$ such that $\phi|_{A\backslash
% \hat{A}} \equiv \psi|_{A\backslash \hat{A}}$.
% \end{lemma}
% \begin{proof}
% MODIFICATIONS NECESSARY
% ??? TODO.
% \end{proof}

% REV IM
\begin{lemma}[Left-cancellation law for $\oplus$]
Let $A,B,C$ be {\sc st}-flow graphs.
\begin{eqnarray*}
A \oplus B = A \oplus C & \Rightarrow & B = C.
\end{eqnarray*}
\label{left-cancellation-law-oplus}
\index{left-cancellation law!for $\oplus$}
\index{$\oplus$!left-cancellation law}
\end{lemma}
\begin{proof}
% REV IM
Let $\langle {A} \oplus {B} \rangle = \langle {A} \rangle
\langle {B} \rangle$ be the canonical $\oplus$-decomposition of
${A} \oplus {B}$, and $\langle {A} \oplus {C} \rangle
= \langle {A} \rangle \langle {C} \rangle$ be the canonical
$\oplus$-decomposition of ${A} \oplus {C}$ respectively.

Fix an isomorphism $\phi$ between $A \oplus B \rightarrow A \oplus C$.
By Proposition
\ref{component-wise-decomposition-of-isomorphisms-under-plus}, $\phi$
maps the elements of $\langle {A} \rangle \langle {B} \rangle$
componentwise to the elements of $\langle {A} \rangle \langle
{C} \rangle$.  It follows that $\phi$ maps the elements of
$\langle {B} \rangle$ componentwise to the elements of $\langle
{C} \rangle$.  Hence a suitable restriction of $\phi$ maps $B$
injectively onto $C$, proving the claim.
\end{proof}

The next lemma is proved in a manner analogous to
Lemma~\ref{left-cancellation-law-oplus}.

\begin{lemma}[Right-cancellation law for $\oplus$]
\label{right-cancellation-law-oplus}
\index{right-cancellation law!for $\oplus$}
\index{$\oplus$!right-cancellation law}
Let $A,B,C$ be flow graphs.
\begin{eqnarray*}
B \oplus A = C \oplus A & \Rightarrow & B = C.
\end{eqnarray*}
\end{lemma}

\begin{proposition}[Commutativity condition for $\oplus$]
\label{plus-commutativity-condition}
\index{$\oplus$!commutativity condition}
% REV IM
Given {\sc st}-flow graphs $A=(G_A, s_A, t_A)$ and $B=(G_B, s_B, t_B)$, 
\begin{eqnarray*}
     A \oplus B & = & B \oplus A
\end{eqnarray*}
iff there exists a flow graph $C$ and integers $k_1$, $k_2$ in
${\mathbb N}$ such that
\begin{eqnarray*}
     A & = & k_1 C, \text{ and}\\
     B & = & k_2 C.
\end{eqnarray*}
\end{proposition}
\begin{proof}
$[\Leftarrow]$ If $A = k_1 C$ and $B = k_2 C$, then $A \oplus B =
(k_1 + k_2) C = B \oplus A$.

$[\Rightarrow]$ The proof is carried by induction on
$\max(|\chi(A)|,|\chi(B)|)$.
Consider the canonical decompositions of $A$ and $B$,
\begin{eqnarray*}
A & = & A^{(0)} \oplus A^{(1)} \oplus A^{(2)} \oplus \cdots \oplus
A^{(|\chi(A)|-1)} \oplus A^{(|\chi(A)|)},\\
B & = & B^{(0)} \oplus B^{(1)} \oplus B^{(2)} \oplus \cdots \oplus
B^{(|\chi(B)|-1)} \oplus B^{(|\chi(B)|)}.
\end{eqnarray*}
If $|\chi(A)|=|\chi(B)|$ then
Proposition~\ref{component-wise-decomposition-of-isomorphisms-under-plus}
tells us that an isomorphism $\phi:A \oplus B \rightarrow B \oplus A$
restricts on the first summand $A$ to yield an isomorphism from $A$ to
$B$.  So in this case, we can take $C=A=B$ and $k_1=k_2=1$.  This
proves the case $\max(|\chi(A)|,|\chi(B)|) = 0$, which forms the basis
of the induction.

Suppose that $\max(|\chi(A)|,|\chi(B)|) > 0$, and $|\chi(A)| \neq
|\chi(B)|$.  Without loss of generality, suppose $|\chi(A)| <
|\chi(B)|$. Then $A^{(i)}=B^{(i)}$ for $i=0, \ldots |\chi(A)|$. It
follows that
\begin{eqnarray}
B^{i + (|\chi(A)| + 1)} & = & B^{(i)} \text{ for } i=0, \ldots,
|\chi(B)| - (|\chi(A)| + 1),\label{exp1}\\
B^{i - (|\chi(A)| + 1)} & = & B^{(i)} \text{ for } i=(|\chi(A)| + 1),
\ldots, |\chi(B)|.\label{exp2}
\end{eqnarray}
If $(|\chi(B)| + 1)$ is divisible by $(|\chi(A)| + 1)$, then
expressions~(\ref{exp1}) and~(\ref{exp2}) above are in fact
equivalent, and in this setting, we take $C=A$, $k_1 = 1$ and $k_2
= \frac{(|\chi(B)| + 1)}{(|\chi(A)| + 1)}$ in order to satisfy the
proposition.
Suppose now that $(|\chi(B)| + 1)$ is {\em not} divisible by
$(|\chi(A)| + 1)$.  Put
\begin{eqnarray*}
d & = & \left\lfloor\frac{(|\chi(B)| + 1)}{(|\chi(A)| +
1)}\right\rfloor, \text{ and }\\
r & = & (|\chi(B)| + 1) \mod (|\chi(A)| + 1).
\end{eqnarray*}
and define
\begin{eqnarray*}
X & = & B^{(0)} \oplus B^{(1)} \cdots \oplus B^{(r-1)},\\
Y & = & B^{(r)} \oplus B^{(r+1)} \cdots \oplus B^{(|\chi(A)|)}.
\end{eqnarray*}
Note that $B = d A \oplus X$ and
\begin{eqnarray*}
         A & = & A^{(0)} \oplus A^{(1)} \oplus A^{(2)} \oplus \cdots
\oplus A^{(|\chi(A)|-1)} \oplus A^{(|\chi(A)|)} \\
           & = & B^{(0)} \oplus B^{(1)} \oplus \cdots B^{(r-1)} \oplus
B^{(r)} \oplus B^{(r+1)} \cdots \oplus B^{(|\chi(A)|-1)} \oplus
B^{(|\chi(A)|)} \\
           & = & X \oplus Y.
\end{eqnarray*}
It follows that $B = d (X \oplus Y) \oplus X$. On the other hand, $X
\oplus Y = A = Y \oplus X$ (see Figure~\ref{xplusy}), since 
\begin{eqnarray*}
X \oplus Y & = & A \\
           & = & A^{(0)} \oplus A^{(1)} \oplus \cdots \oplus
A^{(|\chi(A)|)}\\
           & = & B^{(0)} \oplus B^{(1)} \oplus \cdots \oplus B^{(r)}
\oplus \cdots \oplus B^{(|\chi(A)|-1)} \oplus B^{(|\chi(A)|)}\\
           & = & B^{(r)} \oplus \cdots \oplus B^{(|\chi(B)|-1)} \oplus
B^{(|\chi(B)|)} \oplus B^{(0)} \oplus \cdots \oplus B^{(r-1)} \\
           & = & Y \oplus X.
\end{eqnarray*}
\begin{figure}[htb]
\centering{\mbox{\psfig{figure=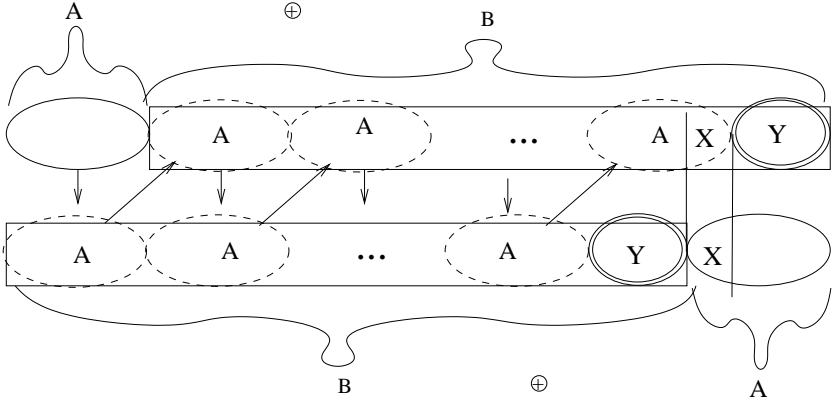}}}
\caption{Inductive step showing $X \oplus Y = A = Y \oplus X$.}
\label{xplusy}
\end{figure}

Since $r \neq 0$ the inductive hypothesis applies to the flow
graphs $X$, $Y$, i.e. there exists some flow graph $Z$ and
suitable integers $l_1$, $l_2$ so that $X = l_1 Z$, $Y = l_2 Z$.
It follows that
\begin{eqnarray*}
A & = & X \oplus Y = l_1 Z \oplus l_2 Z =  (l_1 + l_2)Z \text{ and }\\
B & = & d A \oplus X = d (X \oplus Y) \oplus X = d(l_1 + l_2)Z
\oplus l_1 Z = {((d+1)l_1 + l_2)}Z.
\end{eqnarray*}
So taking $C=Z$, $k_1 = l_1 + l_2$ and $k_2 = (d+1)l_1 + l_2$, the
Proposition is proved.
\end{proof}

\begin{figure}[htb]
\centering{\mbox{\psfig{figure=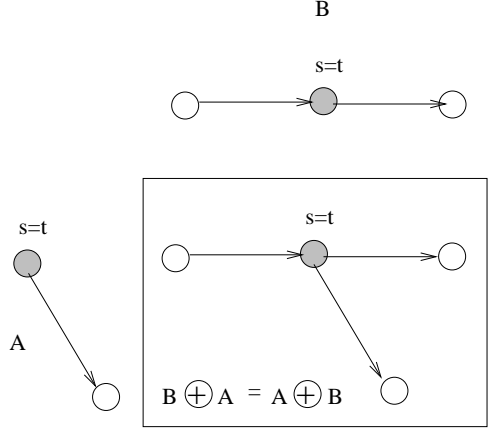}}}
\caption{An example illustrating that the Commutativity Condition fails for general flow graphs.}
\label{comm-cont-counterexample}
\end{figure}

% REV IM
Note that Proposition \ref{plus-commutativity-condition} fails if $A$
and $B$ are both not {\sc st}-flow graphs.  To see this, consider the
example shown in Figure \ref{comm-cont-counterexample}.

%------------------------------------------------------
\subsection{Multiplicative Properties}

In this section we present properties of $\otimes$.
\label{multiplication-section}

\begin{lemma}
\label{edge-otimes-bijection}
Given flow graphs $A$ and $B$, there is a natural bijective
correspondence
\begin{eqnarray*}
  \Lambda_{A, B} : E[A \otimes B] \rightarrow E[G_A] \times E[G_B]
\end{eqnarray*}
\end{lemma}
\begin{proof}
Fix an edge $e$ in $E[G_A \otimes G_B]$.  Then $e$ appears in $(A
\otimes_\eta B)$ at some stage $i$ where $1 \leqslant i \leqslant
|E[G_A]|$ (where $\eta$ is the enumeration specified in
Definition~\ref{otimes-definition}).  We define $\lambda_A(e)$ to be the
edge
$e_i\in E[G_A]$.  At stage $i$ we effectively replace edge
$e_i=(u_i,v_i)$ with a new disjoint copy of $G_B$--by gluing $s_B$
with $u_i$ and $t_B$ with $v_i$.  Thus the edge $e$ corresponds to
some edge $\lambda_B(e)$ in this new disjoint copy of $G_B$.  The
desired bijection $e \mapsto \Lambda_{A, B}(e) = (\lambda_A(e),
\lambda_B(e))$ is thus obtained.  Note that the bijection
$\Lambda_{A, B}$ is independent of the enumeration $\eta$ of the
edges of $E[G_A]$ which appears in the definition of $A \otimes
B$.
\end{proof}

\begin{lemma}[Associativity of $\otimes$]
\label{otimes-associative} \index{$\otimes$!associativity} The
operation $\otimes$ is associative.
\end{lemma}
\begin{proof}
Given flow graphs $A = (G_A, s_A, t_A)$, $B = (G_B, s_B, t_B)$, $C =
(G_C, s_C, t_C)$, we want to show:
\begin{eqnarray*}
(A \otimes B) \otimes C & = & A \otimes (B \otimes C).
\end{eqnarray*}
By Lemma~\ref{edge-otimes-bijection}, the map $\Lambda_{A \otimes B, C}$
is a bijective correspondence between the edges of $(A \otimes B)
\otimes C$ and $(E[G_A] \times
E[G_B]) \times E[G_C]$.  Likewise,
the edges of $A \otimes (B \otimes C)$ are in bijective
correspondence with $E[G_A] \times (E[G_B] \times E[G_C])$, via
$\Lambda_{A, B \otimes C}$.  Obviously $(E[G_A] \times E[G_B])
\times E[G_C]$ is in bijective correspondence with $E[G_A] \times
(E[G_B] \times E[G_C])$ by the map $\pi:((e_1, e_2), e_3) \mapsto
(e_1, (e_2, e_3))$.  Then the composite map
\begin{eqnarray*}
  \mu = \Lambda_{A, B \otimes C} \circ \pi \circ \Lambda_{A \otimes B,
C}
\end{eqnarray*}
is an isomorphism of flow graphs which carries $(A
\otimes B) \otimes C$ to $A \otimes ( B \otimes C)$.
\end{proof}

\begin{example}
\label{non-comm-times-example}
Let $A$ be the flow graph consisting of a directed cycle of length $3$
and let source and target vertices be any two vertices on this cycle.
Then it is easy to check that there is no flow graph isomorphism
between $A\otimes F_2$ and $F_2\otimes A$ (see Figure
\ref{non-comm-times-figure}).  
\begin{figure}[htb]
\centering{\mbox{\psfig{figure=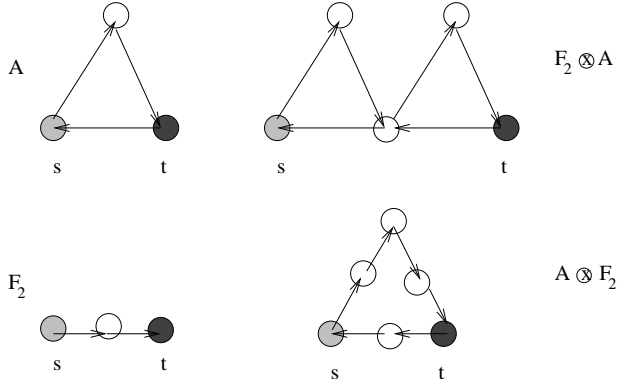}}}
\caption{Example showing the non-commutativity of multiplication in ${\mathcal F}$.}
\label{non-comm-times-figure}
\end{figure}
\end{example}

The previous example proves the next lemma.
\begin{lemma}
\index{$\otimes$!non-commutativity}
The operation $\otimes$ is not commutative.
\end{lemma}

\begin{lemma}[Right-distributivity of $\otimes$ over $\oplus$]
\label{distributive-right}
\index{$\otimes$!right-distributivity over $\oplus$}
For any flow graphs $A, B, C$,
\begin{eqnarray*}
     (A \oplus B) \otimes C & = & (A \otimes C) \oplus (B \otimes C)
\end{eqnarray*}
\end{lemma}
\begin{proof}
Fix $e \in E[(A \oplus B) \otimes C]$.  Then define $\beta_0(e) =
\Lambda_{A \tinyoplus B, C}(e)$.  Note that $\beta_0(e)=(e', f)$, where
$e'$ is an edge in $E[G_{A \tinyoplus B}]$ and $f$ is one in $E[G_C]$.
Define $\beta_1: E[G_{A \tinyoplus B}] \rightarrow E[G_{A}] \cup E[G_{B}]$
so that
\begin{eqnarray*}
\beta_1(e) & = &
\left\{
\begin{array}{rcl}
\sigma^{\NStinyoplus-1}_{t_A \approx s_B}(e) & \text{ if } & e\in
Im(\sigma^{\NStinyoplus-1}_{t_A \approx s_B})\\
\tau^{\NStinyoplus-1}_{t_A \approx s_B}(e) & \text{ if } & e\in
Im(\tau^{\NStinyoplus-1}_{t_A \approx s_B}).
\end{array}
\right.
\end{eqnarray*}
Then $\beta_1 \circ \beta_0$ maps $E[(A \oplus B)
\otimes C]$ injectively into $(E[G_{A}]\times E[G_C]) \cup
(E[G_{B}]\times E[G_C])$.  Define $\beta_2$ by taking
\begin{eqnarray*}
\beta_2(e) & = &
\left\{
\begin{array}{rcl}
\Lambda_{A,C}^{-1}(e) & \text{ if } & e\in E[G_{A\otimes C}]\\
\Lambda_{B,C}^{-1}(e) & \text{ if } & e\in E[G_{B\otimes C}].
\end{array}
\right.
\end{eqnarray*}
Then $\beta_2$ maps $(E[G_{A}]\times E[G_C]) \cup (E[G_{B}]\times
E[G_C])$ into $E[G_{A\otimes C}] \cup E[G_{B\otimes C}]$ injectively.
Finally, define $\beta_3$ by taking
\begin{eqnarray*}
\beta_3(e) & = &
\left\{
\begin{array}{rcl}
\sigma^{\NStinyoplus}_{t_{\!\!A\!\!\otimes\!\! C}\!\!\approx\;
s_{\!\!B\!\!\!\otimes\!\! C}}(e) & \text{ if } & e\in E[G_{
A\!\!\otimes\! C}]\\
\tau^{\NStinyoplus}_{t_{\!\!A\!\!\otimes\!\! C}\!\!\approx\;
s_{\!\!B\!\!\!\otimes\!\!C}}(e) & \text{ if } & e\in E[G_{
B\!\!\otimes\! C}].
\end{array}
\right.
\end{eqnarray*}
Then $\beta_3$ maps $E[G_{A\otimes C}] \cup E[G_{B\otimes C}]$
injectively into
$E[G_{(A\!\!\otimes\!\!C) \tinyoplus (B\!\!\otimes\!\!C)}]$.
The composite map $\beta_3 \circ \beta_2 \circ \beta_1 \circ
\beta_0$ maps the edges of $(A \oplus B) \otimes C$ injectively
into the edges of $(A\otimes C) \oplus (B\otimes C)$, and is the
desired flow graph isomorphism demonstrating the claimed equality.
\end{proof}

Let $A$ be the flow graph consisting of a directed cycle of length $3$
taking source and target vertices to be any two vertices on this
cycle.  Observe that $A \otimes (F_1 \oplus F_1) = A \otimes F_2$,
while $(A \otimes F_1) \oplus (A \otimes F_1) = A \oplus A = 2 A = F_2
\otimes A$.  Referring to Figure \ref{non-comm-times-figure} again, we
see that $A \otimes F_2 \neq F_2 \otimes A$.  Thus, we have shown

\begin{lemma}[Non Left-distributivity of $\otimes$ over $\oplus$]
\label{not-distributive-left}
\index{$\otimes$!not left-distributivity over $\oplus$}
There exist flow graphs $A, B, C$,
\begin{eqnarray*}
     A \otimes (B \oplus C) & \not= & (A \otimes B) \oplus (A \otimes C)
\end{eqnarray*}
\end{lemma}

\begin{definition}
Given flow graphs $A, B$ at least one of which is
non-trivial, and a flow graph $C$,
we say
\begin{eqnarray*}
     A/B=C & \text{  iff  } & A = C\otimes B\\
     A\backslash B=C & \text{  iff  } & A = B\otimes C.
\end{eqnarray*}
If there is no $C$ for which $A/B=C$, we say that $A/B$ does not exist
and $A$ is {\em not right-divisible} by $B$.  If there is no $C$ for
which $A\backslash B=C$, we say that $A\backslash B$ does not exist,
and $A$ is {\em not left-divisible} by $B$.  By convention, we say
that $\ozero / \ozero$ and $\ozero \backslash \ozero$ are undefined.
\end{definition}

Clearly if $m$ and $n$ are standard integers then $F_m/F_n$ iff
$F_m\backslash F_n$ iff $m$ is divisible by $n$.

\begin{lemma}
For all flow graphs $A, B, C$
\begin{eqnarray*}
     A/C & = & (A/B) \otimes (B/C)\\
     A\backslash C & = &  (B\backslash C) \otimes (A\backslash B)
\end{eqnarray*}
whenever these graphs exist.
\end{lemma}
\begin{proof}
Suppose $A/C = K_1$ and $B/C = K_2$.  By definition, $A = K_1
\otimes B$ and $B = K_2 \otimes C$.  Thus, $A = K_1 \otimes (K_2
\otimes C)$, which by Lemma~\ref{otimes-associative} is $(K_1
\otimes K_2) \otimes C$.  Thus $A/C$ exists and equals $K_1
\otimes K_2 = (A/B) \otimes (B/C)$.
Suppose $A\backslash C = K_1$ and $B\backslash C = K_2$.  Then by
definition, $A = B \otimes K_1$ and $B = C \otimes K_2$.  Thus, $A =
(C \otimes K_2) \otimes K_1$, which by Lemma~\ref{otimes-associative}
is $C \otimes (K_2 \otimes K_1)$.  Thus $A\backslash C$
exists and equals $K_2 \otimes K_1 = (B\backslash C) \otimes
(A\backslash B)$.
\end{proof}

\begin{lemma}[Distributivity of right-divisibility over $\oplus$]
For all flow graphs $A, B, C$,
\begin{eqnarray*}
    A/B \oplus C/B & = & (A \oplus C)/B
\end{eqnarray*}
\end{lemma}
\begin{proof}
Suppose $A/B = K_1$ and $C/B = K_2$.  Then by definition, $A = K_1
\otimes B$ and $C = K_2 \otimes B$.  Thus $A \oplus C = (K_1 \otimes B) \oplus (K_2
\otimes B)$ which by Lemma~\ref{distributive-right},
equals $(K_1 \oplus K_2) \otimes B$.  It follows that $(A
\oplus C)/B$ equals $K_1 \oplus K_2$, which is $A/B \oplus C/B$.
\end{proof}

\begin{observation}[Non-distributivity of left-divisibility over $\oplus$] 
Note that Lemma~\ref{not-distributive-left} can be used to construct
examples that demonstrate non-distributivity of left-divisibility over
$\oplus$.  For example, let $B$ be a directed cycle of length $3$ with
any two vertices as $s_B$ and $t_B$.  Take $A = B \otimes F_2$.  Then
$A \backslash B = F_2$.  Now take $C = B$. Then $C\backslash B = F_1$
and so $(A\backslash B) \oplus (C \backslash B) = F_2 \oplus F_1 =
F_3$. Since $A \oplus C \neq B \otimes F_3$, we see that $(A \oplus C)
\backslash B \neq (A\backslash B) \oplus (C \backslash B)$.
\end{observation}

In Definition~\ref{splitting vertex for a flow graph}, we introduced
the notion of a splitting vertex.  Now, in unravelling information
about $\otimes$, we require the notion of a splitting edge.

\begin{definition}[Splitting edge for a flow graph]
\label{splitting edge for a flow graph}
\index{splitting!edge for a flow graph}
\index{edge!splitting --- for a flow graph}
\index{flow graph!edge splitting for a ---}
Let $A = (G_A, s_A, t_A)$ be a flow graph.  A {\em splitting edge} of
$A$
is an edge $e \in E[G_A]$ with the property that $G\backslash e$ has
precisely two components, one of which contains $s_A$ and the other
contains $t_A$.  We denote the set of all splitting edges in $A$ as
$\Delta(A)$.
\end{definition}

\begin{lemma}
\label{oplus-irreducible-delta}
\index{splitting! edges and $\oplus$-irreducibles}
For any flow graph $A$, if $A$ is $\oplus$-irreducible and $A \neq
\oone$, then $\Delta(A)=\emptyset$.
\end{lemma}
\begin{proof}
If $A \neq F_1$ and $e=(u,v)$ is a splitting edge then either $u$ or $v$
or both must be a splitting vertex.  Hence $A$ is $\oplus$-reducible.
\end{proof}

\begin{lemma}[Splitting edges in $\oplus$-decompositions]
\label{splitting edge in oplus decompositions}
Given a flow graph $A$, let $\langle A \rangle$ be its
$\oplus$-decomposition.  Then there exists a map $i:\Delta(A) \rightarrow \{0, 1,
\ldots, \chi(A)\}$ which injectively associates to every splitting
edge an $\oplus$-irreducible component in the $\oplus$-decomposition
of $A$, such that $A^{i(e)}$ is a component consisting only of edge
$e$, and is isomorphic to $F_1$.
\end{lemma}
\begin{proof}
Appealing to Proposition~\ref{decomposition-is-valid}, fix $\phi$ an
isomorphism from $A$ to the component decomposition of $\langle A
\rangle$.  By Proposition~\ref{decomposition-is-valid}, every
$A^{(i_e)}$ is $\oplus$-irreducible.  By
Lemma~\ref{oplus-irreducible-delta} it either has no splitting edges or
it is $F_1$.  Suppose $e$ is a splitting edge in $A$.  Then $\phi(e)$ is
a splitting edge inside $A^{(i_e)}$ for some $i_e$ in $\{0, \ldots,
|\chi(A)|\}$.  It follows that $A^{(i_e)} = F_1$ and the map $i:e
\mapsto i_e$ enjoys the property claimed in the lemma.
\end{proof}

\begin{remark}
\label{cut vertex in product}
Given flow graphs $A$ and $B$, a splitting vertex in $A \otimes B$
comes either from a splitting vertex of $A$ or from a splitting vertex
in (a copy of) $B$ which lies on a splitting edge of $A$.
\end{remark}

\begin{observation}[$\oplus$-decompositions for products]
\label{canonical-plus-decomp-for-product} \index{canonical
$\oplus$-decomposition!for products} Given flow graphs $A, B$,
\begin{eqnarray*}
     \langle A \otimes B \rangle & = & ( A^{(0)}\oplus A^{(1)} \oplus
\cdots \oplus A^{(\chi(A))} )\otimes B \\
& = & (A^{(0)}\otimes B) \oplus (A^{(1)}\otimes B) \oplus \cdots \oplus
(A^{(\chi(A))}\otimes B)\\
& = & \langle A^{(0)}\otimes B \rangle \langle A^{(1)}\otimes B \rangle
\cdots \langle A^{(\chi(A))}\otimes B \rangle.
\end{eqnarray*}
\noindent Since $A^{(i)}$ is $\oplus$-irreducible, by
Lemma~\ref{oplus Irreducibility Lemma} we know that
$\chi(A^{(i)})=\emptyset$.
\noindent Case (i).  $A^{(i)} \neq F_1$.  Then by
Lemma~\ref{decomposition-is-valid}, $\Delta(A^{(i)})=\emptyset$.
It follows from Remark~\ref{cut vertex in product} that
$A^{(i)}\otimes B$ has no splitting vertices, so by
Lemma~\ref{oplus Irreducibility Lemma}, it is
$\oplus$-irreducible. Hence, $\langle A^{(i)}\otimes B \rangle$
is a one-element sequence consisting of $A^{(i)}\otimes B$.
\noindent
Case (ii).  $A^{(i)} = F_1$.  Then $A^{(i)}\otimes B =
B$, so $A^{(i)}\otimes B$ is $\oplus$-irreducible iff $B$ is
$\oplus$-irreducible.  In this case $\langle A^{(i)}\otimes B
\rangle = \langle B \rangle$ is a $(|\chi(B)|+1)$-element
subsequence consisting of the $\oplus$-decomposition of $B$.
\end{observation}

\begin{lemma}[$\oplus$-decomposition length for products]
\label{canonical-plus-decomp-length-for-product}
\index{canonical $\oplus$-decomposition!length for products}
Given flow graphs $A, B$,
\begin{eqnarray*}
     |\chi(A \otimes B)| & = & |\chi(A)| + |\Delta(A)|\cdot|\chi(B)|.
\end{eqnarray*}
\end{lemma}
\begin{proof}
Consider each $\oplus$-irreducible component $A^{(i)}$ in
$\langle A \rangle$.  If $A^{(i)} = F_1$ then by
Observation~\ref{canonical-plus-decomp-for-product}, it contributes
$(|\chi(B)|+1)$ components in $\langle A \oplus B\rangle$. There
are $|\Delta(A)|$ components in $\langle A \rangle$ which are
isomorphic to $F_1$, so these together account for $|\Delta(A)|
\cdot (|\chi(B)|+1)$ components in $\langle A \otimes B
\rangle$.  The remaining $|\chi(A)| + 1 - |\Delta(A)|$ components in
$\langle A \rangle$ (again by
Observation~\ref{canonical-plus-decomp-for-product}) each contribute $1$
component to $\langle A \otimes B \rangle$.  It follows that the
total number of components in $\langle A \otimes B \rangle$ is
\begin{eqnarray*}
|\chi(A)| + 1 - |\Delta(A) |+ |\Delta(A)| (|\chi(B)| + 1) & = & \\
|\chi(A)| + |\Delta(A)| |\chi(B)| + 1.
\end{eqnarray*}
This shows that $|\chi(A \otimes B)| = |\chi(A)| + |\Delta(A)| |\chi(B)|$ as
desired.
\end{proof}

\begin{proposition}[$\oplus$-irreducibility for products]
\label{plus-irreducibilty-for-product}
\index{$\oplus$-irreducibility!for products} Let $A \neq \oone$ be a
flow graph.  Then, $A$ is $\oplus$-irreducible iff for all flow graphs
$B$, $A\otimes B$ is $\oplus$-irreducible.
\end{proposition}
\begin{proof}
$[\Leftarrow]$ Taking $B = \oone$, we see that $A$ is
$\oplus$-irreducible.

$[\Rightarrow]$ If $A$ is $\oplus$-irreducible, then $\chi(A) =
\emptyset$.  Since $A \neq \NSoone$, it follows that $\Delta(A) =
\emptyset$.  So by
Lemma~\ref{canonical-plus-decomp-length-for-product}, for any flow
graph $B$, $\chi(A\otimes B) = \emptyset$.  It follows that $A\otimes
B$ is $\oplus$-irreducible.
\end{proof}

\begin{definition}
\index{prime flow graph}
\index{flow graph!prime}
A flow graph $A$ is called {\em right-prime} if $A/B$ exists only for
$B=\oone$ or $B=A$.  Similarly a flow graph $A$ is called {\em
left-prime} if $A\backslash B$ exists only for $B=\oone$ or $B=A$.
\end{definition}

Note that a natural number $n$ is prime iff the flow graph $F_n$ is
prime.

\begin{lemma}
For all flow graphs $A$, $A$ is right-prime iff $A$ is left-prime.
\end{lemma}
\begin{proof}
Suppose $A$ is right-prime. Assume that $A \backslash B$ exists. We
want to show that $B = F_1$ or $B = A$. Since $A \backslash B$ exists,
it follows $A = B \otimes K$ for some $K$. Thus $A/K$ exists, and
since $A$ is right-prime it follows that $K = F_1$ or $K = A$. If
$K=A$ then $B = F_1$.  If $K = F_1$, then $B = A$.
\end{proof}

%------------------------------------------------------
\subsection{Order Properties}
\label{order-section}

In this section we explore the relationship between strong ordering by
$\stronglesseq$ and weak ordering by $\weaklesseq$.  While the two
orders coincide on the graphical natural numbers, neither order is
anti-symmetric on all of ${\mathcal F}$, and only $\stronglesseq$ is
transitive.  On the other hand, many of the laws that govern the
relationship between $\leqslant$, $+$ and $\times$ in ${\mathcal N}$
continue to hold for $\weaklesseq$, $\oplus$ and $\otimes$ in
${\mathcal F}$, but these laws are violated under the ordering
$\stronglesseq$.

\begin{lemma}[Strong Order Preservation]
\index{strong order!preservation} For flow graphs $A, B, C$, if $A
\stronglesseq B$ then
\begin{eqnarray*}
&& (A\otimes C) \stronglesseq (B\otimes C).
\end{eqnarray*}
\end{lemma}
\begin{proof}
Let $A=(G_A, s_A, t_A)$, $B=(G_B, s_B, t_B)$.  Since $A \stronglesseq
B$ there are graph embeddings $\phi_s:G_A \rightarrow G_B$ and
$\phi_t:G_A \rightarrow G_B$ which satisfy $\phi_s(s_A) = s_B$ and
$\phi_t(t_A) = t_B$.  Define $\gamma_s: E[G_A] \times E[G_C] \rightarrow E[G_B] \times E[G_C]$ by
\begin{eqnarray*}
(e,f) & \mapsto &  (\phi_s(e),f).
\end{eqnarray*}
Then the composite map
\begin{eqnarray*}
\Phi_s^C : E[A\otimes C]
\begin{array}{c}
\Lambda_{A, C}\\
\longrightarrow
\end{array}
E[G_A] \times E[G_C]
\begin{array}{c}
\gamma_s\\
\longrightarrow
\end{array}
E[G_B] \times E[G_C]
\begin{array}{c}
\Lambda_{B, C}^{-1}\\
\longrightarrow
\end{array}
E[B \otimes C]
\end{eqnarray*}
defines an embedding of $G_{A\otimes C} \rightarrow G_{B\otimes C}$
which takes $s_{A\otimes C}$ to $s_{B\otimes C}$.
An analogous construction can be carried out to produce a map
$\Phi_t^C$ which embeds $G_{A\otimes C} \rightarrow G_{B\otimes C}$
and sends $t_{A\otimes C}$ to $t_{B\otimes C}$.
\end{proof}

% REV IM
\begin{lemma}[Strong Order Violations]
\index{strong order!violation}
\label{strong order violation}
There exist {\sc st}-flow graphs $A$, $B$ and $C$ for which
which
\begin{eqnarray*}
(i) (A\oplus C) \notstronglesseq (B\oplus C) \\
(ii) (C\oplus A) \notstronglesseq (C\oplus B) \\
(iii) (C\otimes A) \notstronglesseq (C\otimes B).
\end{eqnarray*}
\end{lemma}
\begin{proof}
See Figure~\ref{strong order violations 1,2,3}.
\end{proof}

\begin{figure}[htb]
\centering{\mbox{\psfig{figure=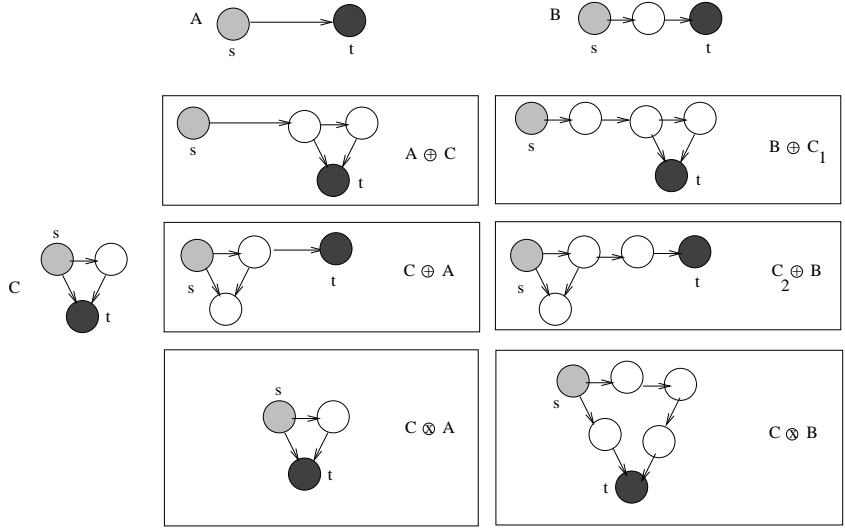}}}
\caption{Strong order violations: (i). $(A\oplus C_1) \notstronglesseq
(B\oplus C_1)$, (ii). $(C_2\oplus A) \notstronglesseq (C_2\oplus
B)$, and (iii). $(C_3\otimes A) \notstronglesseq (C_3\otimes B)$.}
\label{strong order violations 1,2,3}
\end{figure}

We consider possible anti-symmetry of $\stronglesseq$. Suppose $A
\stronglesseq B$ and $B \stronglesseq A$.  There is a graph embedding
$\phi_s:G_A \rightarrow G_B$ which satisfies $\phi_s(s_A) = s_B$.
Hence $|V[G_A]| = |V[\phi_s(G_A)]| \leqslant |V[G_B]|$ and $|E[G_A]| =
|E[\phi_s(G_A)]| \leqslant |E[G_B]|$.  Since $B \stronglesseq A$,
there is a graph embedding $\psi_s:G_B \rightarrow G_A$ which
satisfies $\psi_s(s_B) = s_A$.  So $|V[G_B]| = |V[\psi_s(G_B)]|
\leqslant |V[G_A]|$ and $|E[G_B]| = |E[\psi_s(G_B)]| \leqslant
|E[G_A]|$.  It follows that $\phi_s$ is actually an isomorphism from
$G_A$ to $G_B$ satisfying $\phi_s(s_A) = s_B$.  A similar argument
shows that there is an isomorphism $\phi_t$ from $G_A$ to $G_B$
satisfying $\phi_t(t_A) = t_B$.  To conclude that $A=B$ requires a
single flow graph isomorphism $\pi$ from $A$ to $B$, satisfying {\em
both} $\pi(s_A) = s_B$ and $\pi(t_A) = t_B$.  Indeed in some cases, no
such isomorphism may exist.

\begin{example}
\label{strong-order-not-antisymm-example}
Let $G_A$ be a directed cycle of length $4$, and take $s_A, t_A$ to be
any two vertices in $V[G_A]$ that are distance $2$ apart.  Put $G_B$
isomorphic to $G_A$, taking $s_B, t_B$ to be two vertices in $V[G_B]$
that are distance $1$ apart.  Then it is easy to verify that
$(G_A,s_A,t_A)=A \stronglesseq B=(G_B,s_B,t_B)$ and $B \stronglesseq
A$.  Clearly, however, $A\neq B$ as flow graphs (see
Figure~\ref{non-antisymmetry strong figure}).
\begin{figure}[htb]
\centering{\mbox{\psfig{figure=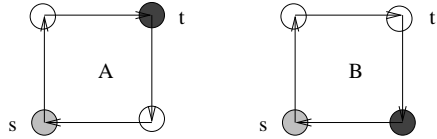}}}
\caption{An example which demonstrates that the strong order is not
antisymmetric.}
\label{non-antisymmetry strong figure}
\end{figure}
\end{example}

The previous example proves the next lemma.

\begin{lemma}[Non-antisymmetry of strong order $\stronglesseq$]
\label{non-antisymmetry strong}
\index{strong order!non-antisymmetry}
\index{antisymmetry!strong order does not satisfy}
There exist flow graphs $A$ and $B$ for which
\begin{eqnarray*}
A \stronglesseq B \text{ and } B \stronglesseq A \text{ but } A \neq B.
\end{eqnarray*}
\end{lemma}

% REV IM
Indeed, since the flow graphs $A$ and $B$ in Example
\ref{strong-order-not-antisymm-example} are both {\sc st}-flow graphs,
it is apparent that antisymmetry of strong order $\stronglesseq$ fails
even for {\sc st}-flow graphs.

\begin{lemma}[Transitivity of strong order $\stronglesseq$]
\index{strong order!transitivity}
For all flow graphs $A, B, C$
\begin{eqnarray*}
A \stronglesseq B \text{ and } B \stronglesseq C \text{ implies } A
\stronglesseq C.
\end{eqnarray*}
\end{lemma}
\begin{proof}
$A \stronglesseq B$: i.e. there are graph embeddings $\phi_s:G_A
\rightarrow G_B$ and $\phi_t:G_A \rightarrow G_B$ which satisfy
$\phi_s(s_A) = s_B$ and $\phi_t(t_A) = t_B$.
$B \stronglesseq C$: i.e. there are graph embeddings $\theta_s:G_B
\rightarrow G_C$ and $\theta_t:G_B \rightarrow G_C$ which satisfy
$\theta_s(s_B) = s_C$ and $\theta_t(t_B) = t_C$.
We want to show $A \stronglesseq C$: i.e. there are graph embeddings
$\alpha_s:G_A \rightarrow G_C$ and $\alpha_t:G_A \rightarrow G_C$
which satisfy $\alpha_s(s_A) = s_C$ and $\alpha_t(t_A) = t_C$.
Put $\alpha_s = \theta_s \circ \phi_s$ and $\alpha_t = \theta_t \circ
\phi_t$.
\end{proof}

\begin{lemma}[Weak Order Preservation]
\index{weak order!preservation} For flow graphs $A, B, C$, if $A
\weaklesseq B$ then
\begin{eqnarray*}
(i) & \text{   } & (A\oplus C) \weaklesseq (B\oplus C) \\
(ii) & \text{   } & (C\oplus A) \weaklesseq (C\oplus B) \\
(iii) & \text{   } & (A\otimes C) \weaklesseq (B\otimes C).
\end{eqnarray*}
\end{lemma}
\begin{proof}
Let $A = (G_A, s_A, t_A), B = (G_B, s_B, t_B)$ and $C = (G_C, s_C,
t_C)$ be given.
$A \weaklesseq B$ implies that there exists an $(s_A,t_A)$-splitting $(H_1, H_2)$
of $G_A$ and graph embeddings
\begin{eqnarray*}
\phi_1 & : & H_1 \rightarrow G_B\\
\phi_2 & : & H_2 \rightarrow G_B
\end{eqnarray*}
satisfy $\phi_1(s_A) = s_B$ and $\phi_2(t_A) = t_B$ and
$\phi_1(E[H_1]) \cap \phi_2(E[H_2]) = \emptyset$.

\vspace{0.1in} 

(i).  Put $K_1 = H_1$ and define $K_2$ to be the
graph obtained by gluing $H_2$ and $G_C$ such that $t_A$ is
identified with $s_C$. Now define $\Phi_1 = \phi_1: H_1
\rightarrow G_B \oplus_{t_B\approx s_C} G_C$, $\Phi_2|_{H_2} =
\phi_2 :H_2 \rightarrow G_B \oplus_{t_B\approx s_C} G_C$ and
$\Phi_2|_{G_C}: G_C \rightarrow G_B \NSoplus_{t_B\approx s_C} G_C$.
Then $(K_1, K_2)$ is an $(s_A, t_C)$-splitting of $G_A
\oplus_{t_A\approx s_C} G_C$ and
\begin{eqnarray*}
\Phi_1 & : & K_1 \rightarrow G_B \oplus_{t_B\approx s_C} G_C,\\
\Phi_2 & : & K_2 \rightarrow G_B \oplus_{t_B\approx s_C} G_C
\end{eqnarray*}
are graph embeddings satisfying $\Phi_1
(s_A) = s_B$ and $\Phi_2 (t_C) = t_C$.

\vspace{0.1in} 

(ii). We define$ L_1$ to be the graph obtained by
gluing $G_C$ and $H_1$ such that $t_C$ is identified with $s_A$
and we put $L_2 = H_2$.
\begin{eqnarray*}
&& \theta_1|_{G_C} : G_C \rightarrow G_C \NSoplus_{t_C\approx s_B}
G_B\\
&& \theta_1|_{H_1} : H_1 \rightarrow G_C \NSoplus_{t_C\approx s_B}
G_B\\
&& \theta_2 = \phi_2:H_2 \rightarrow G_C \NSoplus_{t_C\approx s_B}
G_B.
\end{eqnarray*}
Then $(L_1, L_2)$ is an $(s_C, t_A)$-splitting of $G_C
\oplus_{t_C\approx s_A} G_A$ and
\begin{eqnarray*}
\theta_1& : & L_1 \rightarrow G_C \oplus_{t_C\approx s_B} G_B,\\
\theta_2& : & L_2 \rightarrow G_C \oplus_{t_C\approx s_B} G_B
\end{eqnarray*}
are graph embeddings satisfying $\theta_1 (s_C) = s_C$ and
$\theta_2(t_A)=t_B$.

\vspace{0.1in} 

(iii).  Put
\begin{eqnarray*}
M_1 & = & H_1 \otimes C \\
M_2 & = & H_2 \otimes C.
\end{eqnarray*}
Since $E[H_1] \cup E[H_2] = E[G_A]$ and $E[H_1]) \cap E[H_2] =
\emptyset$, it follows that $E[H_1 \otimes C] \cup E[H_2 \otimes C] =
E[G_A \otimes C]$ and $E[H_1 \otimes C] \cap E[H_2 \otimes C] =
\emptyset$.  Thus $M_1, M_2$ are an $(s,t)$-splitting of $G_A \otimes
C$.
Now take $\beta_1 : E[M_1] \rightarrow E[G_B \otimes G_C]$ and $\beta_2
:
E[M_2] \rightarrow E[G_B \otimes G_C]$ defined by
\begin{eqnarray*}
\beta_1 & : & (e,f)  \mapsto (\phi_1 e, f)\\
\beta_2 & : & (e',f) \mapsto (\phi_2 e', f)
\end{eqnarray*}
for $e \in E[H_1]$, $e' \in E[H_2]$ and $f \in C$.  The
injectivity of $\beta_1$ and $\beta_2$ follows immediately from
injectivity of $\phi_1$ and $\phi_2$.  Since $\phi_1(E[H_1]) \cap
\phi_2(E[H_2]) = \emptyset$, it follows that $\beta_1(E[H_1
\otimes C]) \cap \beta_2(E[H_2 \otimes C]) = \emptyset$.  Since
$\phi_1(s_A) = s_B$ and $\phi_2(t_A) = t_B$, it follows that
$\beta_1(s_{A \otimes C}) = s_{B \otimes C}$ and $\beta_2(t_{A
\otimes C}) = t_{B \otimes C}$. Thus, the maps $\beta_1$ and
$\beta_2$ demonstrate $A \otimes C \weaklesseq B \otimes C$.
\end{proof}

\begin{lemma}[Weak Order Violations]
\index{weak order!violation}
\label{weak order violation}
There exist flow graphs $A, B, C$ for which $A \weaklesseq B$ but
\begin{eqnarray*}
(C\otimes A) \notweaklesseq (C\otimes B).
\end{eqnarray*}
\end{lemma}
\begin{proof}
\begin{figure}[htb]
\centering{\mbox{\psfig{figure=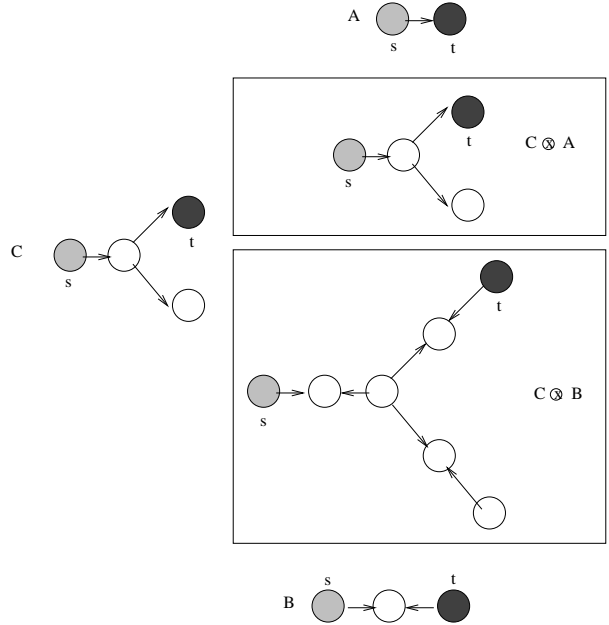}}}
\caption{An example which demonstrates weak order violation: $(C\otimes
A) \notweaklesseq (C\otimes B)$.}
\label{weak order violation figure}
\end{figure}
See Figure~\ref{weak order violation figure}.
\end{proof}

\begin{lemma}[Non-transitive weak order $\weaklesseq$]
\index{weak order!non-transitivity}
There exist flow graphs $A, B, C$
\begin{eqnarray*}
A \weaklesseq B \text{ and } B \weaklesseq C \text{ but } A \notweaklesseq C.
\end{eqnarray*}
\end{lemma}
\begin{proof}
\begin{figure}[htb]
\centering{\mbox{\psfig{figure=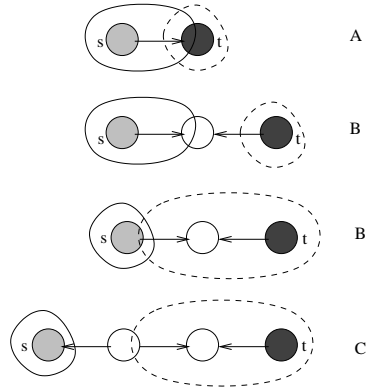}}}
\caption{An example which demonstrates that the weak order is not
transitive.}
\label{non-transitivity weak figure}
\end{figure}
See Figure~\ref{non-transitivity weak figure}.
\end{proof}

% REV IM
\begin{lemma}[Non-antisymmetry of weak order $\weaklesseq$]
\index{weak order!non-antisymmetry}
\index{antisymmetry!weak order does not satisfy}
There exist flow graphs $A$ and $B$ for which
\begin{eqnarray*}
B \weaklesseq A \text{ and } A \weaklesseq B \text{ but } A \neq B 
\end{eqnarray*}
\end{lemma}
\begin{proof}
Since strong order implies weak order, Lemma \ref{non-antisymmetry strong} 
and the example in Figure~\ref{non-antisymmetry strong figure} immediately yield:
\end{proof}

% REV IM
Indeed, since the flow graphs $A$ and $B$ in Example
\ref{strong-order-not-antisymm-example} are both {\sc st}-flow graphs,
it is apparent that antisymmetry of weak order $\stronglesseq$ fails
even for {\sc st}-flow graphs.

%--------------------------------------------------------------------
%--------------------------------------------------------------------
%--------------------------------------------------------------------
\newpage
\section{Conclusions and Future Work}

As we have seen, strikingly many theorems that are true in ${\mathcal
N}$ continue to hold in ${\mathcal F}$, though some fail.  Our future
research program will proceed on two tracks.  

\vspace{0.1in}
Informally, for each ``classical'' theorem $\phi$ in $TA\backslash
Th({\mathcal F})$:
\begin{itemize}
\item[(1)] We shall consider the structure of maximal subsets
$X_\phi$ which have the property that the submodel ${\mathcal X_\phi}
\define ({\mathcal F}|_{X_\phi}) \models \phi$.  Of particular
interest are sets $X_\phi$ which properly contain $i({\mathbb N})$.
\item[(2)] We shall describe a corresponding theorem $\phi'$ in
$Th({\mathcal F})$, such that $\phi' \equiv \phi$ when restricted to
$i({\mathbb N})$.
\end{itemize}

\vspace{0.1in}
{\bf Examples of specific questions include:}

\begin{enumerate}
\item[i.] Characterize flow graph pairs for which antisymmetry of strong
order holds.
\item[ii.] Characterize $\otimes$-commuting pairs, i.e. under what
conditions on flow graphs $A$ and $B$ does $A \otimes B = B
\otimes A$?
\item[iii.] Does $A\otimes B = A\otimes C$ imply $B=C$?  Does $B\otimes
A = C\otimes A$ imply $B=C$?  In other words, does $\otimes$
satisfy a left/right cancellation law?
\item[iv.] {\em Graph Prime Factorization Conjecture.}  Every flow
graph is uniquely expressible (up to some well-defined reordering) as the
product of prime flow graphs.
\item[v.] Describe solution sets (in ${\mathcal F}$) for one-variable
equations having the form $p(x)=q(x)$, where $p$ and $q$ are
polynomials with coefficients from ${\mathcal F}$.
\end{enumerate}

\subsection*{Acknowledgements}

% MOD3
The authors would like to thank Martin Kassabov for his incisive
reading of this paper, insightful comments, and for bringing to our
attention that the Crossing Summands Lemma used in the proof of
Proposition
\ref{component-wise-decomposition-of-isomorphisms-under-plus}) fails
to hold for general flow graphs.  The example demonstrating this (in
Figure \ref{ssl-fails}) is also due to him.

The first author would like to thank the Center for Computational
Science at the Naval Research Laboratory, Washington DC and ITT
Industries for supporting ongoing research efforts in mathematics and
computer science.

\bibliographystyle{siam}
\bibliography{arithX1}

\printindex
\end{document}